\documentclass[a4paper, reqno, 12pt]{amsart}
\usepackage{amsmath, amssymb, eucal, amscd, amstext, enumerate, mathrsfs, amsfonts}
\usepackage[T3,T1]{fontenc}
\DeclareSymbolFont{tipa}{T3}{cmr}{m}{n}
\DeclareMathAccent{\invbreve}{\mathalpha}{tipa}{16}

\newtheorem{thm}{Theorem}[section]
\newtheorem{lem}[thm]{Lemma}
\newtheorem{propdefn}[thm]{Proposition and Definition}
\newtheorem{prop}[thm]{Proposition}
\newtheorem{cor}[thm]{Corollary}

\theoremstyle{definition}
\newtheorem{defn}[thm]{Definition}

\theoremstyle{remark}
\newtheorem{eg}[thm]{Example}
\newtheorem{rem}[thm]{Remark}

\newenvironment{prf}{{\noindent \textbf{Proof:}\ }}{\hfill $\Box$\\ \smallskip}

\numberwithin{equation}{section}

\newcommand{\ti}{\tilde}
\newcommand{\la}{\langle}
\newcommand{\ra}{\rangle}
\newcommand{\smnoind}{\smallskip\noindent}
\newcommand{\CL}{\mathcal{B}}
\newcommand{\CK}{\mathcal{K}}
\newcommand{\CA}{\mathcal{A}}
\newcommand{\CF}{\mathcal{F}}
\newcommand{\proj}{{\rm Proj}^1}
\newcommand{\op}{\operatorname{OP}^1}
\newcommand{\CJ}{\mathcal{J}}
\newcommand{\CS}{\mathcal{SF}}
\newcommand{\D}{\mathcal{D}}
\newcommand{\CT}{\mathcal{II}}

\newcommand{\KI}{\mathfrak{I}}
\newcommand{\KJ}{\mathfrak{J}}
\newcommand{\her}{\operatorname{her}}
\newcommand{\dc}{{\rm d}}
\newcommand{\tw}{{\rm I\!I}}
\newcommand{\ttr}{{\rm I\!I\!I}}
\newcommand{\SF}{{\rm sf}}
\newcommand{\RF}{\mathscr{F}}

\newcommand{\tI}{{\rm I}}
\newcommand{\KH}{\mathfrak{H}}
\newcommand{\al}{{\rm al}}

\newcommand{\pr}{{\rm Prim}}

\newcommand{\T}{\mathscr{T}}
\newcommand{\CI}{\mathcal{I}}
\newcommand{\tsinf}{{\rm ti}}
\newcommand{\hull}{{\rm hull}}
\newcommand{\BC}{\mathbb{C}}
\newcommand{\BN}{\mathbb{N}}

\begin{document}

\baselineskip=1.1\baselineskip
\parskip=0.2\baselineskip
\itemsep=3pt

\title{On the decomposition into Discrete, type $\tw$ and type $\ttr$ $C^*$-algebras}

\author{Chi-Keung Ng  \and Ngai-Ching Wong }

\address[Chi-Keung Ng]{Chern Institute of Mathematics and LPMC, Nankai University, Tianjin 300071, China.}
\email{ckng@nankai.edu.cn}

\address[Ngai-Ching Wong]{Department of Applied Mathematics, National Sun Yat-sen University,
Kaohsiung, 80424, Taiwan.}
\email{wong@math.nsysu.edu.tw}

\thanks{The authors are supported by National Natural Science Foundation of China (11071126 and 11471168), as well as
Taiwan MOST grants (102-2115-M-110-002-MY2 and 104-2115-M-110-009-MY2).}

\date{December 4, 2014; This version: Jun 24, 2016}



\begin{abstract}
We obtained a ``decomposition scheme'' of $C^*$-algebras.
We show that the classes of discrete $C^*$-algebras (as defined by Peligard and Zsid\'{o}),
type $\tw$ $C^*$-algebras  and type $\ttr$ $C^*$-algebras (both defined by Cuntz and Pedersen)
form a good framework to ``classify'' $C^*$-algebras.
In particular, we found that these classes are closed under strong Morita equivalence, hereditary $C^*$-subalgebras as well as taking ``essential extension'' and ``normal quotient''.
Furthermore, there exist the largest discrete finite ideal $A_{\dc,1}$, the largest discrete essentially infinite ideal
$A_{\dc,\infty}$, the largest type $\tw$ finite ideal $A_{\tw,1}$, the largest type $\tw$ essentially infinite ideal $A_{\tw,\infty}$,
and the largest type $\ttr$ ideal $A_\ttr$ of any $C^*$-algebra $A$ such that $A_{\dc,1} + A_{\dc,\infty} + A_{\tw,1} + A_{\tw,\infty} + A_\ttr$ is an essential ideal of $A$.
This ``decomposition'' extends the corresponding one for $W^*$-algebras.

We also give a closer look at $C^*$-algebras with Hausdorff primitive spectrum, $AW^*$-algebras as well as local multiplier algebras of $C^*$-algebras.
We find that these algebras can be decomposed into continuous fields of prime $C^*$-algebras over a locally compact Hausdorff space, with each fiber being non-zero and of one of the five types mentioned above.
%finite discrete, essentially infinite discrete, finite type $\tw$, essentially infinite type $\tw$ or type $\ttr$.
%If, in addition, $A$ is discrete (respectively, essentially infinite), there is an open dense subset of $\Omega$ on which each fiber is discrete %(respectively, essentially infinite).
\end{abstract}
\maketitle

\section{Introduction}

\medskip

Murray and von
Neumann
defined in \cite{Murray-vonNeumann36} (see also \cite{KRII, Murray90}) three types of  $W^*$-algebras according to the abelianness and finiteness properties of their
projections.
Since a $C^*$-algebra needs not have any non-zero projection, a similar
classification for $C^*$-algebras cannot go verbatim.
Cuntz and Pedersen defined (in \cite{CP79}) type $\tw$ and type $\ttr$ $C^*$-algebras according to certain abelianness and finiteness properties of their positive elements.
They used these, together with type $\tI$ $C^*$-algebras, to obtain a classification scheme that captures some features of the $W^*$-algebra counterpart.

\medskip

In  \cite{NW-Mv-class}, we use open projections in $A^{**}$ to obtain another classification scheme for $C^*$-algebras parallel to the one of Murray and von Neumann.
Meanwhile, we also observe that discrete $C^*$-algebras (as defined by Peligard and Zsid\'{o}), type $\tw$ $C^*$-algebras and type $\ttr$ $C^*$-algebras also form a good classification scheme, and some of the results in \cite{NW-Mv-class} have their counterparts in this scheme.
We develop a more comprehensive theory in the current paper.
%The original aim of this article is to presents these results.
%However, we find that one has much better results and a much better classification scheme in this case.
Note that the overlap materials between the current paper and \cite{NW-Mv-class} is not significant.
Actually, only the arguments of Lemma \ref{lem:morita-equiv} and Theorem \ref{thm:decomp} have overlap with the corresponding results in \cite{NW-Mv-class}, and some of the statements that have correspondences in \cite{NW-Mv-class} have a different proofs here.
Moreover, most of the results in this paper are completely new, e.g., all the results in Subsections 3.3 and 3.5, as well as Sections 4 and 6 have no correspondences in \cite{NW-Mv-class} at all.
Conversely, more than half of the results in \cite{NW-Mv-class} has no correspondence in the current paper neither.

\medskip

Let us first recall the notion of discrete, type $\tw$ and type $\ttr$ $C^*$-algebras.
A positive element $x\in A_+$ is said to be \emph{abelian} (in $A$) if $\overline{xAx}$ is an abelian algebra (see \cite[p.191]{Ped79}).
As in the literature, $A$ is said to be \emph{anti-liminary} if there is no non-zero abelian element in $A_+$.
Following Cuntz and Pedersen (\cite{CP79}), for $x,y\in A_+$,  we write $x\sim y$
if there is a sequence $\{z_k\}_{k\in \mathbb{N}}$ in $A$
such that $x =\sum_{k=1}^\infty z_k^*z_k$ and  $y =\sum_{k=1}^\infty z_kz_k^*$ (in norm).
A positive element $x$ is said to be
\emph{finite} in $A$ if one has $y=x$ whenever $y\in A_+$ satisfying $0\leq y\leq x$ and $y\sim x$.
It can be shown that  a projection in a $W^*$-algebra is finite as a projection
in the sense  of Murray-von Neumann if and only if it is finite as a positive element in the sense  of  Cuntz-Pedersen (see e.g. the proof of Proposition \ref{prop:RR0}(d) in Section 3).

\medskip

\begin{defn}[Cuntz-Pedersen \cite{CP79}; Peligard-Zsid\'{o} \cite{PZ00}]
A $C^*$-algebra $A$ is said to be

\smnoind
(a) \emph{discrete}  if every non-zero element in $A_+$ dominates a non-zero abelian element in $A$ (\cite[Definition 2.1]{PZ00}).

\smnoind
(b) \emph{finite}
(respectively, \emph{semi-finite}) if every non-zero element in  $A_+$ is finite in $A$
(respectively, dominates a non-zero positive finite element  in $A$) (\cite[p.140]{CP79}).

\smnoind
(c) \emph{of type $\tw$} if it has no non-zero abelian element and it is semi-finite (\cite[p.149]{CP79}).

\smnoind
(d) \emph{of type $\ttr$} if it has no non-zero finite element (\cite[p.149]{CP79}).
\end{defn}

\medskip

It follows from the definition that any $C^*$-subalgebra of a finite $C^*$-algebra is finite.

\medskip

A $W^*$-algebra $M$ is a type $\tI$, type $\tw$ or type $\ttr$ $W^*$-algebra
if and only if $M$ is a discrete, type $\tw$ or type $\ttr$ $C^*$-algebra (see  Proposition \ref{prop:RR0}(d)).
We also recall that a type $\tI$ $C^*$-algebra is discrete but the converse is not true (e.g.\ $\CL(\ell^2)$).
In fact, we obtained in \cite[Proposition 4.3(a)]{NW-Mv-class}
the following relation:
\begin{quotation}
\emph{A $C^*$-algebra $A$ is of type $\tI$ if and only if $A$ as well as all the primitive quotient $C^*$-algebras of $A$ are discrete.}
\end{quotation}

\medskip

It is not hard to see that $A$ is discrete if and only if every non-zero hereditary
$C^*$-subalgebra of $A$ contains a non-zero abelian hereditary $C^*$-algebra, or equivalently, a non-zero abelian element
(see \cite[Theorem 2.3]{PZ00}).
We obtain alternative looks of discrete, type $\tw$ and type $\ttr$ $C^*$-algebras as in the following theorem (below, by a \emph{normal ideal}, we mean the set of annihilators of another ideal).
Observe that an ideal in a $W^*$-algebra is a normal ideal if and only if it is weak-$^*$-closed.
This theorem summarizes Theorem \ref{thm:sf-C-st-alg}, Remark \ref{rem:main1} and Corollary \ref{cor:type-annih-id}.

\medskip

\begin{thm}\label{thm:main1}  Let $A$ be a $C^*$-algebra.

\smnoind
(a) $A$ is discrete if and only if every non-zero closed ideal (or equivalently, every non-zero normal ideal) of $A$ contains a non-zero abelian hereditary $C^*$-subalgebra.

\smnoind
(b) $A$ is of type $\tw$ if and only if there is no non-zero abelian hereditary
$C^*$-subalgebra of $A$ and every non-zero closed ideal (or equivalently, every non-zero normal ideal) of $A$ contains a non-zero finite hereditary $C^*$-subalgebra.

\smnoind
(c) $A$ is of type $\ttr$ if and only if it contains no non-zero finite hereditary $C^*$-subalgebra.
\end{thm}

\medskip

We will show
that the above three types of $C^*$-algebras are invariance under strong Morita equivalence and under  essential extension.
Here, a $C^*$-algebra $B$ is an \emph{essential extension} of a $C^*$-algebra $A$ if
$A$ is an essential ideal of $B$ (in the sense that $I\cap A\neq (0)$ for every non-zero closed ideal $I$ of $B$).
%The following summaries  Theorem \ref{thm:sf-C-st-alg}(c), Proposition \ref{prop:type-hered}, Theorem \ref{thm:small-class} and Corollary \ref{cor:cp-type-III-anti-fin} in this paper.

\medskip

\begin{thm}\label{thm:main2}
Let $A$ and $B$ be $C^*$-algebras.

\smnoind
(a) Assume that $B$ is either strongly Morita equivalent  to $A$ or
is an essential extension of  $A$.
Then  $A$ is discrete (respectively, of type
$\tw$, of type $\ttr$ or semi-finite) if and only if so is $B$ (see Theorem \ref{thm:sf-C-st-alg}(c) and Proposition \ref{prop:type-hered}(b)).

\smnoind
(b) If $A$ is discrete (respectively, of type
$\tw$, of type $\ttr$ or semi-finite), $B$ is a hereditary $C^*$-subalgebra of $A$, and $J$ is a normal ideal of $A$, then both $B$ and $A/J$ are discrete (respectively, of type
$\tw$, of type $\ttr$ or semi-finite) (see Proposition \ref{prop:type-hered}(b)\&(c)).

\smnoind
(c) The class of discrete (respectively,  semi-finite or  type $\tw$)
$C^*$-algebras is the smallest class that contains all abelian (respectively,  finite or anti-liminary finite)
$C^*$-algebras and is closed under strong Morita equivalence and
essential extension (see Theorem \ref{thm:small-class}).

\smnoind
(d) $A$ is of type $\ttr$ if and only if for every hereditary $C^*$-subalgebra $B\subseteq A$, there is an essential closed ideal of $B$ with the bidual being a properly infinite $W^*$-algebra (see Corollary \ref{cor:cp-type-III-anti-fin}).
\end{thm}

\medskip

We say that a $C^*$-algebra is \emph{essentially infinite} if it does not contain  any non-zero
finite ideal.
By Proposition \ref{prop:equiv-anti-fin}(a), a $C^*$-algebra is essentially infinite if there is an essential closed ideal with its bidual being a properly infinite $W^*$-algebra.
Part (d) above tells us that a $C^*$-algebra is of type $\ttr$ if and only if all of its hereditary $C^*$-subalgebras are essentially infinite.

\medskip

We obtain the following classification scheme of $C^*$-algebras.

\medskip

\begin{thm}\label{thm:main3}
Let $A, B$ be  $C^*$-algebras.

\smnoind
(a) There exists the largest discrete (respectively, semi-finite, type $\tw$ and type $\ttr$)
hereditary $C^*$-subalgebra $A_\dc$ (respectively, $A_\SF$, $A_\tw$ and $A_\ttr$) of $A$ (see Theorem \ref{thm:decomp}(a)).

\smnoind
(b) $A_\dc$, $A_\SF$, $A_\tw$ and $A_\ttr$ are ideals of $A$ such that
$A_\dc$, $A_\tw$ and $A_\ttr$ are disjoint and  $A_\ttr\cap A_\SF = (0)$ (see Theorem \ref{thm:decomp}(a)).

\smnoind
(c) $A_\dc + A_\tw + A_\ttr$ is an essential   ideal of $A$, and $A_\dc + A_\tw$ is an essential   ideal of $A_\SF$ (see Theorem \ref{thm:decomp}(b)).

\smnoind
(d) $A/A_\dc$ (respectively, $A/(A_\tw+A_\ttr)^{\bot\bot}$, $A/A_\ttr$ and $A/A_\SF$) is the universal anti-liminary (respectively, discrete, semi-finite and  type $\ttr$) normal quotient of $A$ (see Corollary \ref{cor:quot}(a) and Theorem \ref{thm:decomp}(c)\&(d)).

\smnoind
(e) If $A$ is semi-finite, then $A/A_\tw$ (respectively, $A/A_\dc$) is the universal discrete (respectively, type $\tw$) normal quotient of $A$ (see Corollary \ref{cor:quot}(b)).
 
\smnoind
(f) If $B$ is a hereditary $C^*$-subalgebra of $A$, then $B_\dc = A_\dc\cap B$, $B_\SF = A_\SF\cap B$, $B_\tw = A_\tw\cap B$ and $B_\ttr = A_\ttr\cap B$ (see Proposition \ref{prop:decomp}(a)).

\smnoind
(g) If $J$ is an essential closed ideal of $A$, then $A_\# = \{x\in A: xJ \subseteq J_\#\}$, for $\#= \dc, \SF, \tw, \ttr$ (see Proposition \ref{prop:decomp}(b)).
 
\smnoind
(h) If $A$ and $B$ are strongly Morita equivalent,
then $A_\#$ and $B_\#$ (for $\# = \dc, \SF, \tw, \ttr$)
corresponds to each other, under the canonical bijection between ideals of $A$
and ideals of $B$ given by an imprimitivity bimodule (see Corollary \ref{cor:morita-equiv-decomp}).

\smnoind
(i) There exist the largest finite ideal $A_1$ and the largest essentially infinite ideal $A_\infty$ of $A$.
One has $A_1\cap A_\infty = (0)$ and $A_1+ A_\infty$ is an essential   ideal of $A$ (see Remark \ref{rem:fin-ess-inf}, Lemma \ref{lem:inf-part} and Theorem \ref{thm:inf-part}(b)).

\smnoind
(j) $A/A_1$ is the universal essentially infinite normal quotient of $A$, and  $A/A_\infty$ is the universal finite normal quotient of $A$ (see Theorem \ref{thm:inf-part}(a)\&(b)).

\smnoind
(k) If $A$ is discrete (respectively, semi-finite or of type $\tw$), then so are $A/A_1$ and $A/A_\infty$ (see Theorem \ref{thm:inf-part}(c)).

\smnoind
(l) Let $J\subseteq A$ be a closed ideal, and $q:A\to A/J^\bot$ be the quotient map.
Then $J_1 = J \cap A_1$, $J_\infty = J\cap A_\infty$,
$q(A)_1 = q(A_1)^{\bot\bot}$ and $q(A)_\infty = q(A_\infty)^{\bot\bot}$ (see Theorem \ref{thm:inf-part}(d)\&(e)).

\smnoind
(m) $A_{\dc,1}:= A_\dc\cap A_1$ (respectively, $A_{\tw,1}:= A_\tw\cap A_1$) is the largest discrete finite (respectively, type $\tw$ finite) ideal of $A$.
On the other hand, $A_{\dc,\infty}:= A_\dc\cap A_\infty$ (respectively, $A_{\tw,\infty}:= A_\tw\cap A_\infty$) is the largest discrete essentially infinite (respectively, type $\tw$ essentially infinite) ideal of $A$.
One has $A_{\dc,1} + A_{\dc,\infty} + A_{\tw,1} + A_{\tw,\infty} + A_\ttr$ being an essential ideal of $A$ (see Corollary \ref{cor:decomp}).

\smnoind
(n) If $A$ is a $W^*$-algebra, then $A_\dc$, $A_\tw$ and $A_\ttr$
are respectively, the type {\tI}, the type $\tw$ and the type $\ttr$ $W^*$-summands of $A$.
Furthermore,  $A_1$ and $A_\infty$ are respectively, the finite part and the infinite part of $A$ (see Theorem \ref{thm:decomp}(e)).

\end{thm}

\medskip

From these, for any $C^*$-algebra $A$, one has
$$A_{\dc,1} \oplus A_{\dc,\infty} \oplus A_{\tw,1} \oplus A_{\tw,\infty} \oplus A_\ttr \subseteq A \subseteq M(A_{\dc,1}) \oplus M(A_{\dc,\infty}) \oplus M(A_{\tw,1}) \oplus M(A_{\tw,\infty}) \oplus M(A_\ttr),$$
and the $C^*$-algebras $M(A_{\dc,1})$, $M(A_{\dc,\infty})$, $M(A_{\tw,1})$, $M(A_{\tw,\infty})$ and $M(A_\ttr)$ are also discrete finite, discrete essentially infinite, type $\tw$ finite, type $\tw$ essentially infinite and type $\ttr$, respectively.

\medskip

As seen in Theorems \ref{thm:main2}(b) and \ref{thm:main3}(c), normal ideals play an important role in the structure theory of $C^*$-algebras.
Hence, prime $C^*$-algebras (i.e., $C^*$-algebra containing no nonzero normal ideal) can be considered as the counterpart of factors in the $C^*$-world.
The following tells us that every prime $C^*$-algebra is of one of the five types as in the above ``decomposition''.

\medskip

\begin{cor}
Any prime $C^*$-algebra is of one of the five types: discrete finite, discrete essentially infinite, type $\tw$ finite, type $\tw$ essentially infinite, or type $\ttr$ (see Proposition \ref{prop:prime>5types}).
\end{cor}

\medskip

\begin{eg}\label{eq:each123}
(a) Discrete finite prime $C^*$-algebras are exactly matrix algebras (see Proposition \ref{prop:prime-discrete}(b)).

\smnoind
(b) Discrete essentially infinite prime $C^*$-algebras are those that contain $\CK(H)$ as an essential ideal for some infinite dimensional Hilbert space $H$ (see Proposition \ref{prop:prime-discrete}(a)).
A concrete example is the unitalization of $\CK(\ell^2)$. 

\smnoind
(c) For a countable ICC group $\Gamma$, its reduced group $C^*$-algebra $C^*_r(\Gamma)$ is a type $\tw$ finite prime $C^*$-algebra (see Example \ref{eg:prime}(b)).

\smnoind
(d) A simple non-type $I$ $AF$-algebra that admits no tracial state is a type $\tw$ essentially infinite prime $C^*$-algebra (see Example \ref{eg:prime}(c)).

\smnoind
(e) The Calkin algebra is a type $\ttr$ prime $C^*$-algebra (see Example \ref{eg:prime}(d)).
\end{eg}

\medskip

One consequence of the above is the following corollary.

\medskip

\begin{cor}
(a)
If  the primitive spectrum $\pr(A)$ of $A$ is Hausdorff, or $A$ is an $AW^*$-algebra, or $A$ is the local multiplier algebra of a $C^*$-algebra, then $A$ can be represented as the algebra of $C_0$-sections of an (F)-Banach bundle over  an open dense subset $\Omega$ of $\pr (ZM(A))$ with each fiber being a nonzero prime $C^*$-algebra of one the five types.
If, in addition, $A$ is discrete (respectively, essentially infinite), there is a dense subset of  $\Omega$ on which each fiber is discrete (respectively, essentially infinite) (see Corollary \ref{cor:ess-simp}).

\smnoind
(b) If the spectrum of $A$ is extremely disconnected (in particular, when $A$ is $AW^*$-algebra), then
$A = A_{\dc,1} \oplus A_{\dc,\infty} \oplus A_{\tw,1} \oplus A_{\tw,\infty} \oplus A_\ttr$ (see Proposition \ref{prop:quasi-Ston}(a)).
\end{cor}

\medskip

\bigskip

\noindent{\it Acknowledgment}

\medskip

The authors would like to thank Lawrence G. Brown, Joachim Cuntz and Huaxin Lin for some discussions and comments on this subject.
They would also like to thank the referee for careful reading of this paper and for comments that helps to improve the presentation of the paper.

\medskip

\section{Notations and preliminaries}

\medskip

Throughout this article, $A$ and $B$ are $C^*$-algebras, $Z(A)$ is the center of $A$, $A_+$ is the positive cone of $A$, $\proj(A)$ is the set of non-zero projections in $A$, $A^{**}$ is the bidual of
$A$ (equipped with the canonical $W^*$-algebra structure),
and $M(A)\subseteq A^{**}$ is the multiplier algebra of $A$.
As usual, we write $ZM(A)$ for $Z(M(A))$.
For any subsets $X,Y,Z\subseteq A$, we write $XY$ (respectively, $XYZ$)
for the linear span of $\{xy:x\in X, y\in Y\}$ (respectively, $\{xyz:x\in X, y\in Y, z\in Z\}$).
%Furthermore, we use $xYz$ %and $XyZ$
%to denote $\{x\}Y\{z\}$
%and $X\{y\}Z$ respectively,
%when $x,z\in A$.
We also write
$$X^\bot\ :=\ \big\{a\in A: aX = \{0\} = Xa\big\}.$$
Notice that if $X$ is an ideal of $A$, then $X^\bot$ is a closed ideal and $X^\bot = \{a\in A: aX = (0)\}$ (note that a closed ideal is $^*$-invariant).

\medskip

A $C^*$-subalgebra  $B$  of $A$ is said to be \emph{hereditary} if
$B_+$ is a \emph{hereditary subcone} of $A_+$ in the following sense:
\begin{equation}\label{eqt:hered-cone}
B_+ = \{a\in A_+: a\leq x, \text{ for some }x\in B_+ \}.
\end{equation}
Clearly, the intersection of two hereditary $C^*$-subalgebras is hereditary (we use the convention that the zero subalgebra is hereditary).
It is well-known that (see e.g. \cite[Theorem 3.2.2]{Murphy}) a $^*$-invariant closed subspace $B\subseteq A$ is a hereditary $C^*$-subalgebra if and only if $BAB\subseteq B$.
Hence, any closed ideal is a hereditary $C^*$-subalgebra. 
Furthermore, if $D\subseteq A$ is any $C^*$-subalgebra, then $D^\bot$ is a hereditary $C^*$-subalgebra.
%It is easy to see that, for every $a\in A$, the  $C^*$-subalgebra $\overline{aAa^*}$ of $A$ is  hereditary.
Moreover, if $D$ is a hereditary $C^*$-subalgebra of $B$ and $B$ is a hereditary $C^*$-subalgebra of $A$,
then $D$ is a hereditary $C^*$-subalgebra of $A$ (in fact, if $a\in A_+$ satisfying $a\leq x\in D_+$, then $a\in B_+$ because $D_+$ is a subset of the hereditary subcone $B_+$ of $A_+$, and hence $a\in D_+$ because $D_+$ is a hereditary subcone of $B_+$).

\medskip

Let us also give a brief account on open projections, which was introduced by Akemann in \cite{Akemann69} (see also
\cite{Akemann70,LinHX90,Ped79,PZ00,Zhang89} for more information).
A projection $p\in A^{**}$ is called an \emph{open projection of $A$} if
there is an increasing net $\{a_i\}_{i\in \mathfrak{I}}$ of positive elements in
$A_+$ with $\lim_i a_i = p$ in the $\sigma(A^{**},A^*)$-topology.
A projection $q\in A^{**}$ is said to be \emph{closed} if $1-q$ is open.
We use $\op(A)$ to denote the collection of non-zero open projections of $A$.
In the case when $A$ is commutative, open projections of $A$ are exactly the images (in $A^{**}$) of characteristic functions of open subsets of the spectrum of $A$.

\medskip

By \cite[Proposition 3.11.9]{Ped79}, we know that a projection $p\in A^{**}$ is open if and only if it lies in the norm-closure of the set $(\ti A_\mathrm{sa})^m$ of $\sigma(A^{**},A^*)$-limits of increasing nets in the self-adjoint part $\ti A_\mathrm{sa}$ of the unitalization $\ti A$ of $A$.
Furthermore, it was shown in \cite[Theorem 3.12.9]{Ped79} that
$$M(A)_\mathrm{sa}  = (\ti A_\mathrm{sa})^m\cap (\ti A_\mathrm{sa})_m,$$
where $M(A)_\mathrm{sa}$ is the self-adjoint part of the multiplier algebra $M(A)$ of $A$ and $(\ti A_\mathrm{sa})_m$ is the set of $\sigma(A^{**},A^*)$-limits of decreasing nets in $\ti A_\mathrm{sa}$.
Thus, any element in $\proj(M(A))\cup \{0\}$ is both an open projection and a closed projection of $A$.
Conversely, if a projection $p\in A^{**}$ is both open and closed, then it is the $\sigma(A^{**},A^*)$-limit of an increasing in $A_+$ as well as a $\sigma(A^{**},A^*)$-limit of a decreasing in $\ti A_\mathrm{sa}$, which means that $p\in M(A)_\mathrm{sa}$.

\medskip

%It is well-known that there is a canonical bijection between open (resp.\ central and open) projections and hereditary $C^*$-subalgebras (resp.\ ideals) for a $C^*$-algebra $A$.
%More precisely,
For every projection (respectively, central projection) $e$ in $A^{**}$ the $C^*$-subalgebra, 
$$
\her_A(e) := eA^{**}e\cap A
$$
is a hereditary $C^*$-subalgebra (respectively, a closed ideal).
Let us list the following two facts from \cite{Ped79} (which may be used implicitly throughout this article):
\begin{enumerate}[O1).]
\item The assignment: $e\mapsto \her_A(e)$ is a bijection from the set of open (respectively, central open) projections on $A$ onto the set of hereditary $C^*$-subalgebras (respectively, closed ideals) of $A$ (see \cite[Remark 3.11.10]{Ped79}).

\item A projection $e\in \proj(A^{**})$ is open in $A$ if and only if $e$ belongs to the  $\sigma(A^{**}, A^*)$-closure of $\her_A(e)$  (see \cite[Proposition 3.11.9]{Ped79}).
In this case, one can choose an increasing net in $\her_A(e)$ that $\sigma(A^{**}, A^*)$-converges to $e$ (see the proof of \cite[Proposition 3.11.9]{Ped79}).
\end{enumerate}

\medskip

The following are some well-known facts about open projections and hereditary $C^*$-subalgebras.

\medskip

\begin{prop}\label{prop:fact-open-hered}
Suppose that $A$ is a $C^*$-algebra and $e,p,q\in \op(A)$ with $p, q\in Z(A^{**})$. %(i.e.\ commutes with every elements in $A^{**}$).

\smnoind
(a) If $z(e)$ is the central cover of $e$ in $A^{**}$, then $z(e)\in \op(A)$ and $\her_A(z(e))$ is the closed ideal, $\overline{A\her_A(e)A}$, generated by $\her_A(e)$.

\smnoind
(b) $ep\in \op(A)$ and $\her_A(ep) = \her_A(e)\cap \her_A(p) = \her_A(e)\her_A(p)\her_A(e)$.

\smnoind
(c) If $\her_A(e)\subseteq \her_A(p) +\her_A(q)$, then  $\her_A(e)\cap \her_A(p) \neq (0)$ or $\her_A(e)\cap \her_A(q) \neq (0)$.
\end{prop}
\begin{prf}
(a) Let $B:= \her_A(e)$ and $I:= \overline{ABA}$.
Denote by $z_0\in \op(A)\cap Z(A^{**})$ the projection with $I = \her_A(z_0) = z_0A^{**} \cap A$.
As $B\subseteq I$, we know that $e\leq z_0$.
If $z\in \proj(A^{**})\cap Z(A^{**})$ satisfying $e\leq z$, then $B\subseteq zA^{**}$ which implies $I\subseteq zA^{**}$ and hence $z_0\leq z$ (see (O2)).
These show that $z_0 = z(e)$.

\smnoind
(b) Suppose that $\{x_i\}_{i\in\KI}$ and $\{y_j\}_{j\in \KJ}$ are increasing nets in $\her_A(e)$ and $\her_A(p)$ that $\sigma(A^{**},A^*)$-converge to $e$ and $p$, respectively.
For a fixed $i\in \KI$, it follows from 
$$ep x_i^{1/2} y_j x_i^{1/2} = x_i^{1/2} y_j x_i^{1/2} = x_i^{1/2} y_j x_i^{1/2}ep$$
(because $p\in Z(A^{**})$) that the increasing net $\big\{x_i^{1/2} y_j x_i^{1/2}\big\}_{j\in \KJ}$ lies  in $\her_A(ep)$, and it will $\sigma(A^{**},A^*)$-converge to $x_ip = x_i^{1/2} p x_i^{1/2}$. 
Consequently, $ep$, being the $\sigma(A^{**},A^*)$-limit of $\{x_ip\}_{i\in \KI}$, lies inside the $\sigma(A^{**},A^*)$-closure of $\her_A(ep)$, and Statement (O2) implies that $ep$ is open in $A$. 

Let $D$ be the hereditary $C^*$-subalgebra $\her_A(e)\cap \her_A(p)$ and $f\in \op(A)$ such that $D = \her_A(f)$.
Since $p$ is central, $ep = e\wedge p$, and we know that $\her_A(ep)\subseteq D$ (and so, $ep\leq f$).
Conversely, as $f\leq e \wedge p$, we know that $D \subseteq \her_A(ep)$.

On the other hand, we denote $D_0 := \her_A(e)\her_A(p)\her_A(e)$.
As every element in $D$ is a product of three elements in $D$, we know that $D\subseteq D_0$.
Conversely, as $\her_A(p)$ is an ideal and $\her_A(e)$ is hereditary, we know that $D_0\subseteq D$.

\smnoind
(c) Since
$
\her_A(p)+\her_A(q)\subseteq \her_A(p+q -pq),
$
we know that $e\leq p+q-pq$.
Suppose on the contrary that
$$
\her_A(e)\cap \her_A(p)\ =\ (0)\ =\ \her_A(e)\cap \her_A(q).
$$
By part (b), we know that $ep = eq = 0$.
Hence, we will arrive at the contradiction that $e = e(p+q-pq) = 0$.
\end{prf}

\medskip

\begin{rem}\label{rem:essetial}
Suppose that $e,f\in \op(A)$ with $e\leq f$ and $B,D\subseteq A$ are hereditary $C^*$-subalgebras.

\smnoind
(a) By \cite[Remark 2.2(b)]{NW-Mv-class}, one may identify
$$\op(\her_A(f)) = \{r\in \op(A): r\leq f\}.$$
Let $\overline{e}^f$ be the closure of $e$ in $\her_A(f)$ (see \cite[Definition II.11]{Akemann69}); i.e., $\overline{e}^f$ is the
smallest closed projection of $\her_A(f)$ that dominate $e$.
Consequently, $f- \overline{e}^f$ is the largest element in
$$\{r\in \op(A)\cup \{0\}: r\leq f \text{ and } re=0\}.$$
From which, we obtain
\begin{equation}\label{eqt:her(e)-bot}
\her_A(e)^\bot = \her_A(1-\overline{e}^1).
\end{equation}
Indeed, one can see from Statement (O2) that $\her_A(e)^\bot = \{x\in A:  ex = 0 = xe\}$.
So, if $p\in \op(A)$, then $\her_A(p)\subseteq \her_A(e)^\bot$ if and only if $pe = 0$, which is equivalent to $p\leq 1-\overline{e}^1$.
Thus, Relation \eqref{eqt:her(e)-bot} follows from the description of $f-\overline{e}^f$ in the above and the fact that $\her_A(e)^\bot$ is a hereditary $C^*$-subalgebra of $A$.

\smnoind
(b) $D$ is said to be \emph{disjoint} from $B$ if $BD = (0)$.
Clearly, $D$ is disjoint from $B$ if and only if $D\subseteq B^\bot$.
As in \cite{Zhang89}, we say that $B$ is \emph{essential} in $A$ if there is no non-zero hereditary $C^*$-subalgebra of $A$ that is disjoint from $B$ (or equivalently, $B^\bot = (0)$).

We say that $e$ is \emph{dense}
in $f$ if $\overline{e}^f = f$.
Relation  \eqref{eqt:her(e)-bot} (applying to the case when $A$ is replaced by $\her_A(f)$) tells us that $e$ is dense in $f$ if and only if $\her_A(e)$ is essential in $\her_A(f)$.

\smnoind
(c) The ideal $I:=\overline{ABA}$ is essential in $A$ as a hereditary $C^*$-subalgebra if and only if for every non-zero closed ideal $J$ of $A$, one has $BJ \neq (0)$ (or equivalently, $B^\bot$ does not contain a non-zero ideal of $A$).
In fact, if $I$ is essential in $A$ but there is a non-zero closed ideal $J\subseteq A$ with $BJ = (0)$, then
$$IJ = \overline{ABA}J \subseteq \overline{ABJ} = (0),$$
which is a contradiction.
Conversely, suppose that $B$ satisfies the said condition, and $D\subseteq A$ is a non-zero hereditary $C^*$-subalgebra.
If $J:= \overline{ADA}$, then $BJ \neq (0)$, and hence 
$$(0)\neq ABADA\subseteq IDA,$$ 
which implies that $ID\neq(0)$ as required.

\smnoind
(d) Let $E$ be a $C^*$-algebra.
Then $A$ contains (a $^*$-isomorphic copy of) $E$ as an essential ideal if and only if there is an injective $^*$-homomorphism $\varphi:A\to M(E)$ such that $E\subseteq \varphi(A)$.

In fact, if $A$ contains $E$ as an essential ideal, then by \cite[Theorem 3.1.8]{Murphy}, there is an injective $^*$-homomorphism from $A$ to $M(E)$ extending the inclusion map $E\subseteq M(E)$.
Conversely, suppose that such a map $\varphi$ exists.
Then $E$ is clearly an ideal of $\varphi(A)$.
Moreover, as $E$ is an essential ideal of $M(E)$ (see e.g. \cite[p.82]{Murphy}), one has
$$\big\{a\in \varphi(A): aE  = \{0\} = Ea \big\} \subseteq \big\{x\in M(E): xE = \{0\} = Ex \big\} = \{0\}.$$
This shows that $E$ is an essential ideal of $\varphi(A)$ (see part (b) above).
\end{rem}

\medskip

Next, we give a brief account for the notion of strong Morita equivalence.
The readers may consult some standard literature on this subject (e.g., \cite{Lan}) for more information.
Let $X$ be a right $A$-module equipped with an ``$A$-valued inner product'', i.e.\ a map $\la \cdot, \cdot\ra_A:X\times X \to A$
such that it is $A$-linear in the second variable and it satisfies
$$\la x, y\ra_A^* = \la y, x\ra_A, \qquad \la x,x \ra_A \geq 0 \qquad (x,y\in X)$$
as well as $\{x\in X: \la x,x \ra_A = 0 \} = \{0\}$.
If $X$ is complete under the norm defined by
$$\|z\|:= \|\la z, z\ra_A\|^{1/2} \qquad (z\in X),$$
then it is called a \emph{Hilbert $A$-module}.
A Hilbert $A$-module is said to be \emph{full} if the linear span of $\{\la x, y\ra_A: x,y\in X\}$ is dense in $A$.
For any $x,y\in X$, we define $\theta_{x,y}: X\to X$ by
$$\theta_{x,y} (z) := x\la y,z\ra_A \qquad (z\in X).$$
The closure, $\CK_A(X)$, of the linear span of $\{\theta_{x,y}:x,y\in X\}$ in the Banach space of bounded linear operators on $X$ is naturally a $C^*$-algebra with involution $\theta_{x,y}^* := \theta_{y,x}$.
If $B$ is another $C^*$-algebra such that there exists a full Hilbert $A$-module $X$ with $B\cong \CK_A(X)$,
then we say that $A$ and $B$ are \emph{strongly Morita equivalent}.

\medskip

\begin{rem}\label{rem:str-Mori-equiv}
Let $X$ be a full Hilbert $A$-module and $B= \CK_A(X)$.

\smnoind
(a) One may equip the conjugate vector space $\ti X$ of $X$ with a Hilbert $B$-module structure:
$$y\ti\ b = (b^*y)\ti\ \quad \text{and} \quad \la x\ti\ , y\ti\ \ra_B := \theta_{x,y} \qquad (x, y\in X; b\in B).$$
For every $a\in A$, the map defined by $\Psi(a)(x\ti\ ) := (xa^*)\ti\ $
belongs to $\CK_B(\ti X)$, and $\Psi$ is a $^*$-isomorphism from $A$ onto $\CK_B(\ti X)$.
In fact, it is easy to see that $\Psi$ is a linear homomorphism from $A$ to the algebra of bounded linear operators on $\ti X$.
As $X$ is full, the relation $\Psi(\la x,y \ra_A) = \theta_{x\ti\ , y\ti\ }$ implies that $\Psi(A)$ is a dense subspace of  $\CK_B(\ti X)$.
Furthermore, $\Psi$ preserves the adjoint since
$$\la \Psi(a^*)(x\ti\ ), y\ti \ \ra = \la (xa)\ti\ ,y\ti\ \ra
= \theta_{xa,y} = \theta_{x,ya^*} = \la x\ti\ , \Psi(a)y\ti \ \ra \qquad
(x,y\in X; a\in A).$$
Hence, $\Psi$ is surjective.
Finally, if $\Psi(a) = 0$, then $a\la x, y \ra a^* = \la xa^*, ya^* \ra = 0$ for any $x,y\in X$ (because $(xa^*)\ti\ = (ya^*)\ti\  = 0$) and the fullness of $X$ implies that $aAa^* = \{0\}$, and hence $a=0$.

\smnoind
(b) One may equip $D:=\Big(\begin{array}{cc}
A & \tilde X\\
X & B
\end{array}\Big)$ with a canonical $C^*$-algebra structure (which is called the \emph{linking algebra}).
%such that
%$$\Big(\begin{array}{cc}
%a & x\ti\ \\
%y & b
%\end{array}\Big)^*\ :=\ \Big(\begin{array}{cc}
%a^* & y\ti\ \\
%x & b^*
%\end{array}\Big)$$
%and
%$$\Big(\begin{array}{cc}
%a & x\ti\ \\
%y & b
%\end{array}\Big)\Big(\begin{array}{cc}
%c & u\ti\ \\
%v & d
%\end{array}\Big)
%\ :=\ \Big(\begin{array}{cc}
%ac + \la x,v\ra_A & \Psi(a)(u\ti\ ) + x\ti\ d\\
%yc + bv & \theta_{y, u} + bd
%\end{array}\Big).$$
Consider $e:= \Big(\begin{array}{cc}
1 & 0\\
0 & 0
\end{array}\Big)\in M(D)$.
Then
$A \cong eDe  = \her_D(e)$ and
$B \cong (1-e)D(1-e) = \her_D(1-e)$.
Moreover, as the closed ideal of $D$ generated by $A$ is
the whole algebra $D$ and the same is true for $B$, we know from Proposition \ref{prop:fact-open-hered}(a) that $z(e) = 1 = z(1-e)$.
As a more explicit reference, the readers may consult, e.g.\ \cite[Theorem II.7.6.9]{Blac}.

\smnoind
(c) For every closed ideal $I$ of $A$, the subset $XI:=\{xa: x\in X; a\in I\} \subseteq X$ is a full Hilbert $I$-module, and
$$\CK_A(XI)\ :=\ \overline{\mathrm{span}\ \!\{\theta_{x,y}:x,y\in XI \}}$$
is an ideal of $B$ that can be identified with
$\CK_I(XI)$.

In fact, if $F$ is the closure of the linear span of $XI$ in $X$, then clearly $F$ is a Hilbert $I$-module, and it follows from the Cohen factorization theorem (see e.g. Theorem 1.11.10 in \cite[p.61]{BD}) that elements in $F$ are of the form $ya$ for some $y\in F\subseteq X$ and $a\in I$, which gives the converse inclusion $F \subseteq XI$.
On the other hand, since any element in $I$ is a product of three elements in $I$, we know from fullness of $X$ that $XI$ is a full Hilbert $I$-module.

It is clear that $\CK_A(XI)$ is an ideal of $B$.
Moreover, one has
$$T(X)\subseteq XI \qquad (T\in \CK_A(XI)),$$
and it is easy to check that $T\mapsto T|_{XI}$ is a surjective $^*$-homomorphism from $\CK_A(XI)$ onto $\CK_I(XI)$.
Suppose that $T|_{XI} = 0$.
If $x\in X$ and $\{u_i\}_{i\in \KI}$ is an approximate unit of $I$, then
$$0 = T(xu_i) = T(x)u_i \to T(x)\in XI,$$
which implies that $T=0$.

Furthermore, part (a) ensures that the map $I\mapsto \CK_I(XI)$ is a bijection from the collection of closed ideals of $A$ onto that of $B$.

\smnoind
(d) Suppose that $C$ is a hereditary $C^*$-subalgebra of $A$ such that $\overline{ACA} = A$.
Then $C$ is strongly Morita equivalent to $A$.
In fact, one may take $X := \overline{CA}$ and equip it with the canonical Hilbert $A$-module structure.
\end{rem}

\medskip

\section{Two news looks of discreteness, type $\tw$ and type $\ttr$}

\medskip

Let us start the main content of this paper with the following simple lemma, which is implicitly included in \cite{CP79}.
Since some of them are not explicitly stated there, we will give a full account here.

\medskip

\begin{lem}\label{lem:finite}
Let $B\subseteq A$ be a non-zero hereditary $C^*$-subalgebra and $\RF^A$ be the set of all non-zero positive elements that are finite in $A$.

\smnoind
(a) If $b\in A$ satisfying $bb^*,b^*b\in B$, then $b\in B$.

\smnoind
(b) $\RF^B = \RF^A\cap B$.

\smnoind
(c) If $x\in \RF^A$ and $\epsilon \in (0, \|x\|)$, then $\overline{(x-\epsilon)_+A(x-\epsilon)_+}$ is a finite $C^*$-algebra.

\smnoind
(d) If $x\in \RF^A$ such that $0$ is an isolated point of $\sigma(x)\cup \{0\}$ (in particular, if $x$ is a projection), then $\overline{xAx}$ is a finite $C^*$-algebra.
\end{lem}
\begin{prf}
	(a) Let $e\in \op(A)$ with $B:=\her_A(e)$.
	As $b^*b\in \her_A(e)$, by considering the polar
	decomposition of $b$, we see that $be = b$.
	Similarly, we have $eb =
	b$.
	Hence, $b= ebe\in \her_A(e)$.
	
	\smnoind
	(b) It is obvious that $\mathscr{F}^A\cap B \subseteq \mathscr{F}^B$.
	Conversely, suppose that $x\in \mathscr{F}^B$.
	Consider $y\in A_+$ and a sequence $\{z_k\}_{k\in \mathbb{N}}$ in $A$ such that $y\leq x$, $y = \sum_{k=1}^\infty z_kz_k^*$ and $x = \sum_{k=1}^\infty z_k^*z_k$.
	Since $B_+$ is a hereditary subcone of $A_+$ (see \eqref{eqt:hered-cone}), we have $y\in B_+$ and $z_k^*z_k, z_kz_k^*\in B_+$ ($k\in \mathbb{N}$).
	By part (a), we know that $z_k\in B$. 
	Now, $x\in \mathscr{F}^B$ gives $y =x$.
	
	\smnoind
	(c) Let $D:= \overline{(x-\epsilon)_+A(x-\epsilon)_+}$ and set
	$$\mathscr{F}_0 := \{a\in A_+\setminus \{0\}: a = ay \text{ for some } y\in \mathscr{F}^A\}.$$
	By \cite[Lemma 4.1]{CP79}, one has $\mathscr{F}_0\subseteq \mathscr{F}^A$.
	Consider $f(t):= \max \{t-\epsilon, 0\}$ ($t\in \sigma(x)$).
	There exists $g\in C(\sigma(x))_+$ as well as $\lambda > 0$ satisfying $g(t)f(t) = f(t)$ and $g(t) \leq \lambda t$  ($t\in \sigma(x)$).
	Hence, $g(x)\leq \lambda x$, and this gives $g(x) \in \mathscr{F}_A$ (because of \cite[Lemma 4.1]{CP79}).
	Moreover, we know from $f(x) = f(x)g(x)$ that for any $z\in D_+$, one has $zg(x) = z$, which implies
	$$z\in \mathscr{F}_0\cap D\subseteq \mathscr{F}^A\cap D \subseteq \mathscr{F}^D.$$
	Consequently, $D$ is a finite $C^*$-algebra.
	
	\smnoind
	(d) As $0$ is an isolated point of $\sigma(x)\cup \{0\}$, one can find $g\in C(\sigma(x)\cup \{0\})_+$ and $\lambda > 0$ such that $g(t)t = t$ and $g(t) \leq \lambda t$  ($t\in \sigma(x)\cup \{0\}$).
	The same argument as in part (c) will give the required conclusion.
\end{prf}

\medskip

The corresponding statement of part (b) above for abelian elements follows from the fact that if $x\in B_+\setminus \{0\}$ is abelian, then $xAx \subseteq x^{1/2}(x^{1/2}Ax^{1/2})x^{1/2} \subseteq \overline{x^{1/2}Bx^{1/2}} = \overline{xBx}$ and hence, is an abelian algebra.

\medskip

Our next lemma is a crucial step toward the new looks and the classification scheme.
A different proof of the statement concerning abelian elements can be found in \cite[Proposition 3.8(a)]{NW-Mv-class}.

\medskip

\begin{lem}\label{lem:morita-equiv}
Let $A$ and $B$ be two strongly Morita equivalent $C^*$-algebras.
Then $A$ has a non-zero abelian element  (respectively, finite elements) if and only if $B$ does.\end{lem}
\begin{prf}
We will only give the proof for the statement concerning finite elements since the other statement
follows from a similar argument.
Let $D$ and $e\in \proj(M(D))\subseteq \op(D)$ be as in Remark \ref{rem:str-Mori-equiv}(b).
In particular, $A = \her_D(e)$ and $z(e) = 1$ in $\op(D)$.

Suppose that $\RF^B\neq \emptyset$.
Lemma \ref{lem:finite}(b) implies that $\RF^D\neq \emptyset$.
By Lemma \ref{lem:finite}(c), we know that $D$ contains a non-zero finite hereditary $C^*$-subalgebra $D_0$.
Let $p\in \op(D)$ with $D_0 = \her_D(p)$.
Then \cite[Theorem 1.9]{PZ00} gives $e_0,e_1,p_0,p_1\in \op(D)\cup \{0\}$ such that $e_0,e_1\leq e$, $p_0, p_1\leq p$, $e_0e_1=0$, $p_0p_1=0$,
$$\overline{e_0 + e_1}^e =e, \ \ \overline{p_0 +
p_1}^p =p, \ \ z(e_0)z(p_0) = 0 \ \text{ and }\ \her_D(e_1)\cong \her_D(p_1)$$
(notice that the equivalence relation in \cite{PZ00} produces a $^*$-isomorphism between $\her_D(e_1)$ and $\her_D(p_1)$).

Let us first show that $p_1\neq 0$.
Suppose on the contrary that $p_1 =0$.
Then $e_1 = 0$ (see Statement (O1)) and $z(e_0)$ is
dense in $z(e) = 1$ (by \cite[Lemma 1.8]{PZ00}).
This implies that $z(p_0) = 0$ (because of Proposition \ref{prop:fact-open-hered}(a) and Remark \ref{rem:essetial}(a)) and hence $p_0=0$, which is absurd because $p_0$ is dense in the non-zero open projection $p$.

Now, $\her_D(p_1)\subseteq D_0$ is a non-zero finite hereditary $C^*$-subalgebra of $D$.
Since $\her_D(e_1)\cong \her_D(p_1)$, we conclude that $\her_D(e_1)$ is a non-zero finite hereditary $C^*$-subalgebra
of $\her_D(e)$.
Consequently, $\RF^A \neq \emptyset$ (by Lemma \ref{lem:finite}(b)) as required.
\end{prf}

\medskip

\subsection{The first set of new looks}

The following theorem is our first set of new looks.
Note that there is no ambiguity in part (b) of this theorem because of Lemma \ref{lem:finite}(b).

\medskip

\begin{thm}\label{thm:sf-C-st-alg}
(a) $A$ is discrete if and only if every non-zero closed ideal of $A$ contains a non-zero abelian element.

\smnoind
(b) $A$ is semi-finite if and only if every non-zero closed ideal of $A$ contains a non-zero finite element.

\smnoind
(c) Let $A$ and $B$ be two strongly Morita equivalent $C^*$-algebras.
If $A$ is discrete (respectively, anti-liminary, semi-finite, of type $\tw$ or of type $\ttr$), then so is $B$.
\end{thm}
\begin{prf}
(a) It is clear that if $A$ is discrete, the required condition holds.
Conversely, suppose that $A$ satisfies the condition as stated, and $D\subseteq A$ is a non-zero hereditary $C^*$-algebra.
If $I := \overline{ADA}$, then Remark \ref{rem:str-Mori-equiv}(d) tells us that $I$ is strongly Morita equivalent to $D$.
Now, we know from Lemma \ref{lem:morita-equiv} that $D$ contains a non-zero abelian element and hence $A$ is discrete (by \cite[Theorem 2.3]{PZ00}).

\smnoind
(b) Since the ``only if part'' is clear (recall that the cone of a closed ideal of $A$ is a hereditary subcone of $A_+$), we only need to consider the ``if part''.
Let $a\in A_+\setminus \{0\}$ and set $D:= \overline{aAa}$.
By Remark \ref{rem:str-Mori-equiv}(d), $D$ is strongly Morita equivalent to $\overline{ADA}$.
The hypothesis and Lemma \ref{lem:morita-equiv} produces $x\in \RF^D = \RF^A\cap D$ with $\|x\|=1$ (see Lemma \ref{lem:finite}(b)).
Notice that if $a^{1/2}x^{1/2} = 0$, then we have a contradiction that $Dx^{1/2} = 0$.
Moreover, since $x^{1/2}ax^{1/2} \leq \|a\| x$ and
$a^{1/2}xa^{1/2} \sim x^{1/2}ax^{1/2}$, we know from \cite[Lemma 4.1]{CP79} that $a^{1/2}xa^{1/2}\in \RF^A$.
As $a^{1/2}xa^{1/2}\leq a$, we see that $A$
is semi-finite.

\smnoind
(c) It follows from Lemma \ref{lem:morita-equiv} that if $A$ is anti-liminary (respectively, type $\ttr$), then so is $B$.
Let $J$ be a non-zero closed ideal of $B$.
By Remark \ref{rem:str-Mori-equiv}(c), there is a closed ideal $I\subseteq A$ that is strongly Morita equivalent to $J$.
If $A$ is discrete (respectively, semi-finite), then by part (a) (respectively, part (b)) as well as Lemma \ref{lem:morita-equiv}, we know that $J$ contains a non-zero abelian (respectively, finite) element, and hence $B$ is discrete (respectively, semi-finite).
Since a $C^*$-algebra is of type $\tw$ if and only if it is anti-liminary and semi-finite, we know that the property of being type $\tw$ is also preserved under strong Morita equivalence.
\end{prf}

\medskip

\begin{rem}\label{rem:main1}
(a) It is clear that one can replace ``abelian element'' in Theorem \ref{thm:sf-C-st-alg}(a) by ``abelian hereditary $C^*$-subalgebra''.
Moreover, by Lemma \ref{lem:finite}(c), one may also replace ``finite element'' in Theorem \ref{thm:sf-C-st-alg}(b) by ``finite hereditary $C^*$-subalgebra''.
The same reason also tells us  that:
\begin{quotation}
a $C^*$-algebra is of type $\ttr$ if and only if it has no non-zero finite hereditary $C^*$-subalgebra.
\end{quotation}

\smnoind
(b) Remark \ref{rem:str-Mori-equiv}(d) and Theorem \ref{thm:sf-C-st-alg}(c) implies that every discrete (respectively, anti-liminary, semi-finite, type $\tw$ or type $\ttr$) hereditary $C^*$-subalgebra of $A$ is contained in a discrete (respectively, anti-liminary, semi-finite, type $\tw$ or type $\ttr$) closed ideal of $A$.
\end{rem}

\medskip

\subsection{Some permanence properties}

\medskip

With the help of Theorem \ref{thm:sf-C-st-alg}, we will obtain
in Proposition \ref{prop:type-hered} below some permanence properties of discreteness,
type $\tw$, type $\ttr$ and semi-finiteness.
Apart from part (a), these properties seem to be new.
Notice that, by Remark \ref{rem:essetial}(c), the hypothesis  in part (b)
means that $B^\bot$ does not contain any non-zero ideal of $A$, but the presentation in the statement here seems more informative.

\medskip

\begin{prop}\label{prop:type-hered}
Let $B$ be a discrete (respectively, type
$\tw$, type $\ttr$, semi-finite or anti-liminary) $C^*$-algebra.

\smnoind
(a) If $A$ is a hereditary $C^*$-subalgebra of $B$, then $A$ is discrete (respectively, of type $\tw$, of type $\ttr$, semi-finite or anti-liminary).

\smnoind
(b) If $B$ is a hereditary $C^*$-subalgebra of a $C^*$-algebra $A$ such that
$I:=\overline{ABA}$ is an essential ideal of $A$, then $A$ is discrete (respectively, of type
$\tw$, of type $\ttr$, semi-finite or anti-liminary).
In particular, $M(B)$ is discrete (respectively, of type $\tw$, of type
$\ttr$, semi-finite or anti-liminary).

\smnoind
(c) If $B$ is a closed ideal of a $C^*$-algebra $A$, then both $B^{\bot\bot}$ and $A/B^\bot$ are discrete
(respectively, of type $\tw$, of type $\ttr$, semi-finite or anti-liminary).

\smnoind
(d) If $I$ is an ideal of $B$, then $B/I^\bot$ is discrete
(respectively, of type $\tw$, of type $\ttr$, semi-finite or anti-liminary).
\end{prop}
\begin{prf}
(a) This part follows from the definitions, Lemma \ref{lem:finite}(b) as well as the paragraph following Lemma \ref{lem:finite}.

\smnoind
(b) We will only consider the case when $B$ is of type $\tw$ (as the argument for the other cases are similar).
By Remark \ref{rem:str-Mori-equiv}(d) and Theorem \ref{thm:sf-C-st-alg}(c), we know that $I$ is of type $\tw$.
Suppose that $(0)\neq D\subseteq A$ is a hereditary $C^*$-subalgebra.
Then $ID$ is non-zero, because $I$ is essential.
Hence, $D\cap I = DID\neq (0)$ and is a hereditary $C^*$-subalgebra of $I$.
This ensures that $D$ contains a non-zero finite element.
On the other hand, if $D = \overline{aAa}$ for an abelian element $a\in A_+$, then $D\cap I$
is a non-zero abelian hereditary $C^*$-subalgebra  of $I$, which is impossible.
These show that $A$ is of type $\tw$.

\smnoind
(c) Let $\varphi: A\to M(B)$ be the canonical $^*$-homomorphism.
As $\ker\varphi = B^\bot$, we know that $\varphi|_{B^{\bot\bot}}$ is injective and so is the induced $^*$-homomorphism $\hat \varphi:A/B^\bot \to M(B)$.
Since both $\varphi(B^{\bot\bot})$ and $\hat \varphi(A/B^\bot)$ contains the image of $B$ as an essential ideal (see Remark \ref{rem:essetial}(d)), the conclusion follows from part (b).

\smnoind
(d) By part (a), we know that $I$ is discrete
(respectively, of type $\tw$, of type $\ttr$, semi-finite or anti-liminary).
Hence, by part (c), we know that $B/I^\bot$ is discrete
(respectively, of type $\tw$, of type $\ttr$, semi-finite or anti-liminary).
\end{prf}

\medskip

\begin{defn}\label{defn:annih-id}
(a) We say that $A$ is an \emph{essential extension} of $B$ if it contains $B$ as an essential ideal.

\smnoind
(b) An ideal $I\subseteq A$ is called a \emph{normal ideal} if $I = J^\bot$ for some ideal $J\subseteq A$.

\smnoind
(c) If $I\subseteq A$ is a normal ideal, then $A/I$ is called an \emph{normal quotient of $A$}.
\end{defn}

\medskip

Notice that the name ``normal ideals'' comes from the fact that an ideal of a $W^*$-algebra $M$ is a normal ideal if and only if it is $\sigma(M,M_*)$-closed.

\medskip

Proposition \ref{prop:type-hered}(b)  implies that the properties of being
discrete, type $\tw$, type $\ttr$, semi-finite and anti-liminary are stable under essential extensions.
Furthermore, Proposition \ref{prop:type-hered}(d) states that these types are preserved under taking normal quotients.
In contrast, the quotient of a discrete $C^*$-algebra by an arbitrary closed ideal needs not be discrete (e.g.\ $\CL(\ell^2)$ is a discrete $C^*$-algebra but $\CL(\ell^2)/\CK(\ell^2)$ is of type $\ttr$).

\medskip

\begin{cor}\label{cor:type-annih-id}
$A$ is discrete (respectively, semi-finite) if and only if every non-zero normal ideal of $A$ contains a non-zero abelian (respectively, finite) hereditary $C^*$-subalgebra.
\end{cor}
\begin{prf}
If $A$ is discrete (respectively, semi-finite), then $A$ will satisfy the said condition because of Theorem \ref{thm:sf-C-st-alg}(a) (respectively, Theorem \ref{thm:sf-C-st-alg}(b)).
Conversely, assume that the said condition is statisfied.
Let $I$ be a non-zero closed ideal of $A$.
The assumption tells us that $I^{\bot\bot}$ contains a non-zero abelian (respectively, finite) element.
By the definition (respectively, Lemma \ref{lem:finite}(c)), $I^{\bot\bot}$ contains a non-zero abelian (respectively, finite) hereditary $C^*$-subalgebra $D$.
As $I$ is an essential ideal of $I^{\bot\bot}$, we know that $I\cap D = DID \neq (0)$.
Now, Remark \ref{rem:main1}(a) implies that $A$ is discrete (respectively, semi-finite).
\end{prf}

\medskip

\subsection{The second set of new looks}

The following can be regarded as another new looks of discreteness, semi-finiteness and type $\tw$.

\medskip

\begin{thm}\label{thm:small-class}
Let $\CA$ (respectively, $\CF$ and $\CF_{al}$) be the class of all abelian (respectively, finite and anti-liminary finite) $C^*$-algebras.
Let $\D$ (respectively, $\CS$ and $\CT$) be the class obtained from $\CA$ (respectively, $\CF$ and $\CF_{al}$) by taking strong Morita equivalences and essential extensions.
Then $\D$ (respectively, $\CS$ and $\CT$) is the class of discrete (respectively, semi-finite and type $\tw$) $C^*$-algebras.
\end{thm}
\begin{prf}
By Theorem \ref{thm:sf-C-st-alg}(c) and Proposition \ref{prop:type-hered}(b),
%(which actually implies Theorem \ref{thm:decomp}(e)),
we know that every member in $\D$ (respectively, $\CS$ and $\CT$) is discrete (respectively, semi-finite and type $\tw$).

Suppose that $A$ is a non-zero discrete $C^*$-algebra.
Then \cite[Theorem 2.3(iii)]{PZ00} gives an abelian hereditary $C^*$-subalgebra $A_1$ of $A$ such that $\overline{AA_1A}$ is an essential ideal of $A$.
Thus, by Remark \ref{rem:str-Mori-equiv}(d), we know that $A\in \D$.

Secondly, assume that $A$ is a semi-finite $C^*$-algebra.
In order to show that $A$ is a member of $\CS$, it suffices to show that there is a finite hereditary $C^*$-subalgebra $A_2\subseteq A$
with  $\overline{AA_2A}$ being an essential ideal of $A$ (because of Remark \ref{rem:str-Mori-equiv}(d)).
%By Lemma \ref{lem:finite}(c), $\CF_A$ is non-empty.
Indeed, by Zorn's Lemma, there exists a maximal family $\{B_i\}_{i\in \KI}$  of  finite hereditary  $C^*$-subalgebras of $A$
such that $\overline{AB_iA}$ and $\overline{AB_jA}$ are disjoint for all $i\neq j$, and we set
$$
A_2\ :=\ \overline{{\sum}_{i\in \KI} B_i}.
$$
Because $A_2$ is the $c_0$-direct sum of finite $C^*$-algebras, it is finite.
Moreover, the ideal $\overline{AA_2A}$ is essential, by the maximality of $\{B_i\}_{i\in \KI}$ and Remark \ref{rem:main1}(a).

Finally, let $A$ be a type $\tw$ $C^*$-algebra.
As in above, we consider a maximal family $\{D_i\}_{i\in \KI}$ of anti-liminary finite hereditary $C^*$-subalgebras with $\overline{AD_iA}$ and $\overline{AD_jA}$ being disjoint when $i\neq j$.
If $A_3 := \overline{{\sum}_{i\in \KI} D_i}$, then $A_3$ is an anti-liminary finite hereditary $C^*$-subalgebra of $A$.
Again by Remark \ref{rem:main1}(a) and the maximality, the ideal $\overline{AA_3A}$ is essential in $A$, and we know that $A$ is a member of $\CT$.
\end{prf}

\medskip

\subsection{$C^*$-algebras of real rank zero and $W^*$-algebras}

\medskip

We want to compare discreteness, type $\tw$ and type $\ttr$ of $C^*$-algebras with the corresponding properties of $W^*$-algebras.
These comparisons could be known, but since we do not find them in the literature, we present them in Proposition \ref{prop:RR0}(d) below for later reference.

\medskip

\begin{prop}\label{prop:RR0}
Let $A$ be a $C^*$-algebra of real rank zero.

\smnoind
(a) $A$ is discrete (respectively, semi-finite) if and only if every non-zero projection in $A$ dominates a non-zero projection that  is an abelian element of $A$ (respectively, is finite in $A$).

\smnoind
(b) $A$ is of type $\tw$ if and only if every non-zero projection in $A$ is not abelian but dominates a
non-zero projection that is finite in $A$.

\smnoind
(c) $A$ is of type $\ttr$ if and only if it has no non-zero projection which is finite in $A$.

\smnoind
(d) If, in addition, $A$ is a $W^*$-algebra, then $A$ is
a type {\tI}, type $\tw$, type $\ttr$ or semi-finite $W^*$-algebra
if and only if it is, respectively, a discrete, type $\tw$, type $\ttr$ or semi-finite $C^*$-algebra.
\end{prop}
\begin{prf}
(a) Suppose that $A$ is discrete (respectively, semi-finite).
Any projection $p\in A$ will dominate a non-zero abelian (respectively, finite) element $x\in A_+$.
As $(pAp)_+$ is a hereditary subcone of $A$, we know that $x\in pAp$ (respectively, $(x- p\|x\|/2 )_+\in pAp$).
Hence, $pAp$ will contains a non-zero abelian (respectively, finite) hereditary $C^*$-subalgebra (by Lemma \ref{lem:finite}(c)), and any non-zero projection in this subalgebra will be dominated by $p$.

Conversely, assume that the said condition holds.
If $I\subseteq A$ is a non-zero closed ideal and $p\in \proj(I)$, then $p$ will dominate an non-zero abelian (respectively, finite) element and this element will belong to $I$ (as $I_+$ is a hereditary subcone of $A_+$).
Thus, $A$ is discrete (respectively, semi-finite) by Theorem \ref{thm:sf-C-st-alg}.

\smnoind
(b) This part follows directly from part (a) (note that if $x\in A_+$ is abelian then any projection in $\overline{xAx}$ is abelian).

\smnoind
(c) The ``only if part'' is clear, and the converse follows from Lemma \ref{lem:finite}(c). 

\smnoind
(d) By parts (a), (b) and (c) above, it suffices to show that a projection is finite in $A$ if and only if it is finite in the sense of Murry and von Neumann.
In fact, one of the implication is clear because the equivalence relation $\sim$ in the beginning of Section 1 is weaker than the Murry-von Neumann equivalence.
Conversely, suppose that a projection $p\in \proj(A)$ is finite in the sense of Murry and von Neumann.
Then $pAp$ is finite as a $W^*$-algebra, and the set of normal tracial states will separate $0$ from other positive elements in $pAp$ (i.e., for any $x\in (pAp)_+\setminus \{0\}$, there exists a normal tracial state $\tau$ with $\tau(x)\neq 0$). 
Hence, we know from \cite[Theorem 3.4]{CP79} that $pAp$ is a finite $C^*$-algebra and $p\in pAp$ is finite in $A$.
\end{prf}

\medskip

\subsection{Purely infinite $C^*$-algebras and tracially infinite $C^*$-algebras}

Another application of Proposition \ref{prop:type-hered} is our next result, which shows that type $\ttr$ is weaker than the notion of purely infiniteness
as defined by Cuntz in \cite{Cu81} (in the case of simple $C^*$-algebras) as well as by Kirchberg and R\o rdam \cite{KR00} (in the general case).
Note that Proposition \ref{prop:pure-inf>III} implies, in particular, \cite[Proposition 4.4]{KR00}.

\medskip

In order to obtain this comparison, we will investigate another property that looks very similar to type $\ttr$ (compare Proposition \ref{prop:pure-inf>III}(a) as well as Corollary \ref{cor:cp-type-III-anti-fin} in the next section), and lies between type $\ttr$ and pure infiniteness.
This property will also be considered in Section 6 below.
Before introducing this property, let us first recall the notion of traces on $C^*$-algebras.

\medskip

As in \cite[\S 5.2]{Ped79}, by a \emph{trace}, we mean a weight $\tau:A_+\to [0,\infty]$ satisfying
$$\tau(u^*xu) =\tau(x)$$
for any $x\in A_+$ and unitary $u$ in the unitalization of $A$.
A trace $\tau$ is said to be \emph{lower semi-continuous} if $\{x\in A_+: \tau(x)\leq \lambda \}$ is norm closed in $A$, for every $\lambda\in \mathbb{R}_+$.
Moreover, as in \cite[Definition 6.1.1]{Dix77}, $\tau$ is called \emph{semi-finite} if
$$\tau(x) = \sup\{\tau(y): 0\leq y\leq x \text{ and } \tau(y)<+\infty\} \qquad (x\in A_+).$$
Note that this semi-finiteness is slightly stronger than the one in \cite{CP79}.
In the following, we denote by $\T(A)$ the set of all lower semi-continuous semi-finite traces on $A$,
and by $T_s(A)$ the set of tracial states on $A$.
%Clearly, the constant null function $0$ belongs to $\T(A)$.

\medskip

\begin{propdefn}\label{prop:pure-inf>III}
$A$ is said to be \emph{tracially infinite} if it contains an essential closed ideal $J$ with $\T(J) = \{0\}$.

\smnoind
(a) $A$ is tracially infinite if and only if $A$ contains an essential ideal $J$ such that for every hereditary $C^*$-subalgebra $B\subseteq J$, it bidual $B^{**}$ is a properly infinite $W^*$-algebra.

\smnoind
(b) If $A$ is tracially infinite, then $A$ is of type $\ttr$.
In particular, every purely infinite $C^*$-algebra is of type $\ttr$.
\end{propdefn}
\begin{prf}
(a) It suffices to show that for any $C^*$-algebra $D$, one has $\T(D) = \{0\}$ if and only if $T_s(B) = \emptyset$ for every non-zero hereditary $C^*$-subalgebra $B\subseteq D$ (we may then apply this statement to $D=J$; note that $B^{**}$ is a properly infinite $W^*$-algebra if and only if $T_s(B) = \emptyset$).

Suppose on the contrary that $\T(D) = \{0\}$ but there exist a non-zero hereditary $C^*$-subalgebra $B\subseteq D$ and an element $\varphi\in T_s(B)$.
Let $\tau$ be the lower semi-continuous trace extension of $\varphi$ as given in the proof of \cite[Lemma 4.6]{CP79}; namely,
$$
\tau(x):= \sup \{\varphi(y): y\in B_+, y\sim z\leq x\} \quad (x\in D_+).
$$
It is not hard to check that $\tau$ is semi-finite, and this gives a contradiction.

Conversely, if $\tau\in \T(D)\setminus \{0\}$, then the semi-finiteness of $\tau$ will produce an element $x\in D_+$ with $0 < \tau(x) < +\infty$.
Without loss of generality, assume that $\|x\| = 1$.
It is clear that $\|x - (x-\frac{1}{n})_+\| \to 0$ and the lower semi-continuity of $\tau$ tells us that $\tau\big((x-\frac{1}{n})_+\big) \to \tau(x)$.
Set $y:=(x-\frac{1}{n_0})_+$, where $n_0$ is big enough so that $\tau\big((x-\frac{1}{n_0})_+\big) > 0$.
Using a similar argument as that of Lemma \ref{lem:finite}(c), one can find $g\in C(\sigma(x))_+$ and $\lambda >0$ such that
$y = yg^{1/2}(x)$ and $g(x) \leq \lambda x$.
Thus,
$$\tau(z)\ =\ \tau\big(g^{1/2}(x)zg^{1/2}(x)\big)\ \leq\ \lambda \|z\| \tau(x)\ <\ +\infty \qquad (z\in (\overline{yDy})_+),$$
and hence $T_s(\overline{yDy}) \neq \emptyset$.

\smnoind
(b) Let $J$ be an essential closed ideal of $A$ with $\T(J) = \{0\}$.
By Proposition \ref{prop:type-hered}(b), it suffices to show that $J$ is of type $\ttr$.
Suppose on the contrary that $J$ contains a non-zero finite hereditary $C^*$-subalgebra $B$ (see Remark \ref{rem:main1}(a)).
Then \cite[Theorem 3.4]{CP79} implies that $T_s(B)\neq \emptyset$, and the argument of part (a) produces a non-zero element in $\T(J)$, which is a contradiction.

For the second statement, if $A$ is purely infinite, it follows from \cite[Proposition 5.1]{KR00} that $\T(A) = \{0\}$ and $A$ is tracially infinite.
\end{prf}

\medskip

Recall that pure infiniteness passes to quotients (see \cite[Theorem 4.19]{KR00}) but the quotient of a type $\ttr$ $C^*$-algebra can be semi-finite (see e.g. \cite[Remark 3.14]{CP79}).
Consequently, these two properties are difference.
On the other hand, we do not know whether tracial infiniteness coincides with type $\ttr$.
One can find in \cite{Ng-ext-lscsf-tr} some equivalences of the statement: ``every type $\ttr$ $C^*$-algebra is tracially infinite''.

\medskip

The following is a direct application of Proposition \ref{prop:pure-inf>III}, Theorem \ref{thm:sf-C-st-alg}(c) and Proposition \ref{prop:type-hered}(b).

\medskip

\begin{cor}\label{cor:pure-inf}
Let $A$ be a purely infinite $C^*$-algebra.

\smnoind
(a) If $A$ is strongly Morita equivalent to a $C^*$-algebra $B$, then $B$ is of type $\ttr$.

\smnoind
(b) An essential extension of $A$ (in particular, $M(A)$) is of type $\ttr$.
\end{cor}

\medskip

Notice that if $A$ and $B$ are strongly Morita equivalent separable $C^*$-algebras  with $A$ being purely infinite, then \cite[Theorem 4.23]{KR00} implies that $B$ is also purely infinite.
We do not know if the same is true for non-separable $C^*$-algebras.

\medskip

\begin{eg}\label{eg:III-not-pi}
For every non-unital $C^*$-algebra $B$, the $C^*$-algebras $M(B\otimes \mathcal{O}_\infty)$ and $M(B\otimes \mathcal{O}_2)$ are of type $\ttr$ (by Corollary \ref{cor:pure-inf}(b) and \cite[Proposition 4.5]{KR00}).
\end{eg}

\medskip

\section{Essentially infinite $C^*$-algebras and another equivalent form of type $\ttr$ $C^*$-algebras}

\medskip

As seen in Theorem \ref{thm:small-class}, stability under both strong Morita equivalence as well as essential extension plays an important role in understanding discrete and type $\tw$ $C^*$-algebras.
Clearly, finiteness is not preserved under strong Morita equivalence (e.g.\ $\mathbb{C}$ is strongly Morita equivalent to $C^*$-algebra $\CK(\ell^2)$ that is not finite).
However, it is stable under essential extension, as stated in part (a) of the following lemma.
This lemma could be known, but we present it here for completeness.

\medskip

\begin{lem}\label{lem:finite-under-extension}
(a) If $B$ is a finite $C^*$-algebra, then every essential extension of $B$ is finite.

\smnoind
(b) If $I$ is a non-zero closed ideal of $A$, every $\tau\in T_s(I)$ extends to an element of $T_s(A)$.
\end{lem}
\begin{prf}
(a) Since a  $C^*$-subalgebra of a finite $C^*$-algebra is finite, we need only to consider the essential extension $M(B)$ of $B$ (see Remark \ref{rem:essetial}(d)).
By \cite[Theorem 3.4]{CP79}, $T_s(B)$ separates $0$ from other positive elements in $B$; in other words,
\begin{equation*}
\{x\in B_+: \tau(x) = 0, \text{ for all }\tau\in T_s(B) \} = \{0\}.
\end{equation*}
For any element $\tau\in T_s(B)$, let $(\pi_\tau, \KH_\tau, \xi_\tau)$ be its GNS construction.
If $\pi_0:= \bigoplus_{\tau\in T_s(B)} \pi_\tau$, then the above tells us that $\pi_0$ is faithful and hence extends to a faithful $^*$-representation $\ti \pi_0$ of $M(B)$.
On the other hand, each $\sigma\in T_s(B)$ extends to a tracial state $\bar \sigma$ on $D := \bigoplus_{\tau\in T_s(B)}^{\ell^\infty} \pi_\tau(B)''$ such that
$$\bar \sigma\big((y_\tau)_{\tau\in T_s(B)}\big)
\ =\ \langle y_\sigma \xi_\sigma, \xi_\sigma\rangle  \qquad \big((y_\tau)_{\tau\in T_s(B)}\in D\big).$$
As $\ti \pi_0(M(B))\subseteq D$,  we can define $\ti \sigma := \bar \sigma\circ \ti \pi_0$.
It is easy to see that $(\pi_{\ti\sigma}, \KH_{\ti \sigma})$ is the canonical extension of $(\pi_\sigma, \KH_\sigma)$ on $M(B)$.
Moreover, we have $\ker \pi_{\ti \sigma}=\{x\in M(B): \ti \sigma(x^*x) = 0\}$, since $\ti \sigma$ is a tracial state.
Now, if $a\in M(B)_+$ satisfying $\ti \sigma(a) = 0$ for every $\sigma\in T_s(B)$, then $\pi_{\ti \sigma}(a^{1/2}) = 0$ ($\sigma\in T_s(B)$), which gives $\ti \pi_0(a^{1/2}) = 0$ and hence $a^{1/2} = 0$.
Consequently, $M(B)$ is finite (by \cite[Theorem 3.4]{CP79}).

\smnoind
(b) As $I$ is an essential closed ideal of $I^{\bot\bot}$, the argument of part (a) tells us that $\tau$ can be extended to a tracial state $\tau_0$ on $I^{\bot\bot}$.
Furthermore, as $I^{\bot\bot}+I^\bot$ is an essential ideal of $A$,
the trivial extension of $\tau_0$ on $I^{\bot\bot}+I^\bot$ (i.e.\ it vanishes on $I^\bot$) will again extends to a tracial state on $A$.
\end{prf}

\medskip

Motivated by the classification theory of  $W^*$-algebras, we make the following definition of essentially infiniteness.
Note that a $W^*$-algebra is essentially infinite if and only if it is properly infinite, but the name ``properly infinite $C^*$-algebras'' is already in used in the literature.
The term ``essentially infinite'' also comes from Proposition \ref{prop:equiv-anti-fin}(a) below. 

\medskip

\begin{defn}
A $C^*$-algebra $A$ is said to be \emph{essentially infinite} if it contains no non-zero finite closed ideal.
\end{defn}

\medskip

\begin{rem}\label{rem:fin-ess-inf}
(a) We set
$$A_\infty\ := \ \bigcap \{J^\bot: J \text{ is a finite closed ideal of }A\} \quad \text{and} \quad A_1:=A_\infty^\bot.$$
It is easy to see that $A_\infty$ is essentially infinite.
On the other hand, if $I$ is an arbitrary essentially infinite ideal of $A$, then for any finite ideal $J\subseteq A$, one has $I\cap J = \{0\}$, which implies that $I\subseteq A_\infty$.

\smnoind
(b) If $I\subseteq A$ is a non-zero closed ideal such that $I\cap A_\infty = \{0\}$, then by the definition, $I \subseteq A_1$.
This shows that $A_1 + A_\infty$ is an essential ideal of $A$.

\smnoind
(c) It is clear that any non-zero closed ideal of an essentially infinite $C^*$-algebra is essentially infinite. 

\smnoind
(d) If $A$ is simple, then either $A = A_\infty$ or $A = A_1$.
\end{rem}

\medskip

We will show that $A_1$ and $A_\infty$ can be viewed as the ``finite part'' and the ``infinite part'' of $A$, respectively.
Let us first give the following lemma.

\medskip

\begin{lem}\label{lem:inf-part}
Let $\CJ$ be a maximal collection of pairwise disjoint non-zero finite ideals of $A$.
Then $J_0:= \overline{\sum_{J\in \CJ} J}$ is a finite ideal of $A$ such that $J_0+ A_\infty$ is essential in $A$.
Moreover, $A_\infty = J_0^\bot$.
\end{lem}
\begin{prf}
$J_0$ is finite because it is a $c_0$-direct sum of finite $C^*$-algebras.
Suppose that $I$ is a non-zero closed ideal of $A$.
If $I$ contains a non-zero finite ideal, then $I\cap J_0\neq (0)$ (otherwise, $\CJ$ cannot be a maximal family).
If $I$ does not contain a non-zero finite ideal, then $I\subseteq A_\infty$ (see Remark \ref{rem:fin-ess-inf}(a)).
This shows that $J_0+A_\infty$ is essential.

Furthermore, it is clear that $A_\infty \subseteq J_0^\bot$.
Conversely, as $J_0^\bot$ cannot contain any non-zero finite ideal (because $\CJ$ is maximal), $J_0^\bot$ is essentially infinite and is contained in $A_\infty$.
\end{prf}

\medskip

Observe that the construction in Lemma \ref{lem:inf-part}
needs not give a maximal finite ideal of $A$.
For example, if $A=\ell^\infty$ and $\CJ$ is the collection of all one-dimensional ideals of $A$, then $J_0 = c_0$,
while the abelian $C^*$-algebra $\ell^\infty$ is itself finite.
%However, we have the following result.

\medskip

\begin{thm}\label{thm:inf-part}
Let $A$ be a $C^*$-algebra, $I \subseteq A$ be a closed ideal of $A$ and $q:A\to A/I^\bot$ be the quotient map.

\smnoind
(a)  $A_1^\bot$ coincides with the largest essentially infinite ideal $A_\infty$ of $A$, and $A/A_1$ is the \emph{universal essentially infinite normal quotient of $A$}; i.e., if a normal quotient of $A$ (see Definition \ref{defn:annih-id}(c)) is essentially infinite, then the corresponding quotient map will factor through $A/A_1$.

\smnoind
(b) $A_1$ is the largest finite ideal of $A$, and $A/A_\infty$ is the \emph{universal finite normal quotient of $A$}; i.e., the quotient map of any finite normal quotient of $A$ factors through $A/A_\infty$.

\smnoind
(c) If $A$ is discrete (respectively, of type $\tw$, semi-finite or anti-liminary), then so are $A/A_1$ and $A/A_\infty$.

\smnoind
(d) $I_1 = I \cap A_1$ and $I_\infty = I \cap A_\infty$.

\smnoind
(e) $q(A_1)\cong A_1/(I^\bot)_1$ (respectively, $q(A_\infty) \cong A_\infty/(I^\bot)_\infty$) and is an essential ideal of $(A/I^\bot)_1$ (respectively, $(A/I^\bot)_\infty$).

\smnoind
(f) If $A$ is a $W^*$-algebra, then $A_1$ and $A_\infty$ are, respectively, the finite part and the infinite part of $A$.
\end{thm}
\begin{prf}
(a) We learn from Remark \ref{rem:fin-ess-inf}(a) that $A_\infty$ is the largest essentially infinite ideal of $A$. 
Suppose that $\theta:A \to M(A_\infty)$ is the canonical embedding.
Then $\ker\theta= A_\infty^\bot =A_1$ and $\theta$ induces a $^*$-monomorphism $\hat \theta:A/A_1\to M(A_\infty)$.
Since $\hat \theta(A/A_1)$ contains $\theta(A_\infty)$ as an essential ideal (see Remark \ref{rem:essetial}(d)), if there exists a non-zero finite ideal $J$ of $A/A_1$, then $\hat\theta(J)\cap \theta(A_\infty)$ will be a non-zero finite ideal of $\theta(A_\infty)$ which is a contradiction.
Thus, $A/A_1$ is essentially infinite.

Assume that $A/I^\bot$ is essentially infinite.
Then $I$ is $^*$-isomorphic to the ideal $q(I)$ of $A/I^\bot$ and hence is essentially infinite (see Remark \ref{rem:fin-ess-inf}(c)).
Therefore, $I \subseteq A_\infty$ and $A_1\subseteq I^\bot$, which means that $q$ factors through $A/A_1$.

On the other hand, $\theta$ restricts to an injection from $A_1^\bot = A_\infty^{\bot\bot}$ to $M(A_\infty)$.
It follows from Remark \ref{rem:essetial}(d) that $\theta(A_\infty)$ is an essential ideal of $\theta(A_1^\bot)$, and the same argument as the essential infiniteness of $A/A_1$ implies that $A_1^\bot$ is essentially infinite.
Consequently, we obtain the converse inclusion $A_1^\bot \subseteq A_\infty$.

\smnoind
(b) Let $\CJ$ and $J_0$ be as in Lemma \ref{lem:inf-part}.
Suppose that $\varphi: A\to M(J_0)$ is the canonical $^*$-homomorphism.
Since $\varphi|_{A_1}$ is a $^*$-monomorphism (because $\ker \varphi = J_0^\bot = A_\infty$), we know that $A_1$ is finite (see Lemma \ref{lem:finite-under-extension}(a) and Remark \ref{rem:essetial}(d)).
Moreover, if $J\subseteq A$ is a finite ideal, then $J\cap A_\infty = \{0\}$ (by the definition of $A_\infty$), and $J\subseteq A_\infty^\bot = A_1$.

To show the finiteness of $A/A_\infty$, we note that the canonical map from $A$ to $M(A_1)$ induces an injection on $A/A_\infty$ (as $A_\infty = A_1^\bot$ by part (a)).
Hence, we know from the finiteness of $A_1$ and Lemma \ref{lem:finite-under-extension}(a) that $A/A_\infty$ is finite.

Assume that $A/I^\bot$ is finite.
Then $q$ restricts to an injection on $I$ and the ideal $q(I)$ of $A/I^\bot$ is again finite.
This shows that $I\subseteq A_1$ and hence $A_\infty =  A_1^\bot \subseteq I^\bot$.
Consequently, $q$ factors through $A/A_\infty$.

\smnoind
(c) By part (a) and the definition of $A_1$, both $A/A_1$ and $A/A_\infty$ are normal quotient of $A$. 
Thus, this part follows from Proposition \ref{prop:type-hered}(d).

\smnoind
(d) If $J$ is any closed ideal of $A$, then, by considering an approximate unit of $I$, one has
\begin{equation}\label{eqt:orth-comp}
J^\bot \cap I = (J\cap I)^\bot \cap I.
\end{equation}
This gives
$$A_\infty \cap I
\ =\  \bigcap \{(J\cap I)^\bot\cap I: J \text{ is a finite ideal of }A\}\\
%& = & \bigcap \{K^\bot\cap I: K \text{ is a finite ideal of }I\}
\ =\ I_\infty,$$
and $A_1\cap I = A_\infty^\bot \cap I = (A_\infty \cap I)^\bot \cap I = I_\infty^\bot \cap I= I_1$.

\smnoind
(e) By part (d) and Relation \eqref{eqt:orth-comp}, one has
$$(I^\bot)_1 = I^\bot \cap A_1 = (I_1)^\bot \cap A_1$$
and the canonical $^*$-homomorphism from $A_1$ to $M(I_1)$ induces a $^*$-monomorphism from
$$q(A_1) \cong A_1/(I^\bot \cap A_1) = A_1/(I^\bot)_1 = A_1/((I_1)^\bot \cap A)$$
to $M(I_1)$.
It follows from Lemma \ref{lem:finite-under-extension}(a) and Remark \ref{rem:essetial}(d) that $A_1/((I_1)^\bot \cap A)$ (whose faithful image in $M(I_1)$ contains $I_1$) is a finite $C^*$-algebra. 
Consequently, $q(A_1)\subseteq q(A)_1$ (because of part (b)).

Suppose that $J$ is a non-zero closed ideal of the finite ideal $q(A)_1\subseteq q(A)$.
Since $q$ restricts to an injection from $q^{-1}(J)\cap I^{\bot\bot}$ to the finite $C^*$-algebra $J$, we know that $q^{-1}(J)\cap I^{\bot\bot}$ is an ideal of $A_1$ and is non-zero (because $q^{-1}(J)\nsubseteq I^\bot$).
Thus, $J\cap q(A_1)\neq (0)$, and $q(A_1)$ is an essential ideal of $q(A)_1$.

In a similar way, the canonical map from $A_\infty$ to $M(I_\infty)$ induces a $^*$-monomorphism from
$$q(A_\infty) \cong A_\infty/(I^\bot\cap A_\infty) = A_\infty/(I^\bot)_\infty = A_\infty/((I_\infty)^\bot\cap A_\infty) $$
(see \eqref{eqt:orth-comp}) to $M(I_\infty)$. 
Remark \ref{rem:essetial}(d) tells us that the image of  $q(A_\infty)$ in $M(I_\infty)$ contains $I_\infty$ as an essential ideal.
Since $I_\infty$ is essentially infinite, so is $q(A_\infty)$ (by the definition) and hence $q(A_\infty)\subseteq q(A)_\infty$ (because of part (a)).

Consider $J'$ to be a non-zero closed ideal of the essentially infinite ideal $q(A)_\infty$ of $q(A)$.
Since $q$ restricts to an injection from $q^{-1}(J')\cap I^{\bot\bot}$ to an ideal of the essentially infinite $C^*$-algebra $J'$ (see Remark \ref{rem:fin-ess-inf}(c)), we know that $q^{-1}(J')\cap I^{\bot\bot}$ is a non-zero (because $q^{-1}(J')\nsubseteq I^\bot$) closed ideal of $A_\infty$.
Thus, $J'\cap q(A_\infty)\neq (0)$, and $q(A_\infty)$ is an essential ideal of $q(A)_\infty$.

\smnoind
(f) Recall that a $W^*$-algebra is finite in the sense of Murray-von Neumann if and only if it is finite in the sense of Cuntz-Pedersen (see the proof of Proposition \ref{prop:RR0}(d)).
As $J^\bot$ is $\sigma(A,A_*)$-closed for every ideal $J\subseteq A$, we see that $A_1$ is the largest $\sigma(A,A_*)$-closed finite ideal of $A$ (as $A_1 = A_\infty^\bot$).
Moreover, $A_\infty$ is the largest $\sigma(A,A_*)$-closed ideal of $A$ that contains no $\sigma(A,A_*)$-closed finite ideal (see part (a)).
\end{prf}

\medskip

\cite[Theorem 3.4]{CP79} tells us that $A$ is finite if and only if it has enough tracial states; in the sense that $\bigcap_{\tau \in T_s(A)} \ker \pi_\tau = (0)$ (we use the convention that this intersection is $A$ when $T_s(A) = \emptyset$).
In the following, we will establish that $A$ is essentially infinite if and only if $\bigcap_{\tau\in T_s(A)} \ker \pi_\tau$ is an essential ideal of $A$.
Notice that it is possible for an essentially infinite $C^*$-algebra to have a tracial state.
For example, if $A$ is the unitalization of $\CK(\ell^2)$, then $A$ is essentially infinite (by Proposition \ref{prop:equiv-anti-fin}(a) below),  but the canonical non-degenerate one-dimensional representation $\tau$ of $A$ is a tracial state.

\medskip

\begin{prop}\label{prop:equiv-anti-fin}
(a) The following statements are equivalent.
\begin{enumerate}
\item $A$ is essentially infinite.
\item $\bigcap_{\tau\in T_s(A)} \ker \pi_\tau$ is an essential ideal of $A$.
\item There is an essential ideal $J$ of $A$ with $J^{**}$ being a properly infinite $W^*$-algebra.
\end{enumerate}

\smnoind
(b) If $A$ is separable, then $A$ is essentially infinite if and only if  $\ker \pi_\tau$ is an essential ideal of $A$ for every $\tau \in T_s(A)$
\end{prop}
\begin{prf}
(a) Set $J_0 := \bigcap_{\tau\in T_s(A)} \ker \pi_\tau$.

\smnoind
$1) \Rightarrow 2)$.
Suppose on the contrary that $J_0^\bot \neq (0)$.
As the representation  $\bigoplus_{\tau\in T_s(A)} \pi_\tau$ is injective on $J_0^\bot$, one knows that
\begin{equation*}\label{eqt:sep-posit-elem}
\{x\in (J_0^\bot)_+: \tau(x) = 0, \text{ for all }\tau\in T_s(A) \} = \{0\}.
\end{equation*}
Thus, we obtain the contradiction that $J_0^\bot$ is a finite closed ideal of $A$ (using \cite[Theorem 3.4]{CP79}).

\smnoind
$2) \Rightarrow 3)$.
We will verify $J_0^{**}$ being a properly infinite $W^*$-algebra.
Assume on the contrary that one can find $\tau_0\in T_s(J_0)$.
By Lemma \ref{lem:finite-under-extension}(b), $\tau_0$ extends to an element $\tau_1\in T_s(A)$.
Since $\pi_{\tau_1}(J_0) = \{0\}$, one has $\tau_0((J_0)_+) = \tau_1((J_0)_+) = \{0\}$, which gives the contradiction that $\tau_0 = 0$.

\smnoind
$3) \Rightarrow 1)$.
Suppose on the contrary that $A$ contains a non-zero finite closed ideal $I$.
Then $I_0:= I \cap J\neq (0)$.
As $I_0$ is finite, it has a tracial state $\tau$ and $\tau$ extends to a tracial state on $J$ (because of Lemma \ref{lem:finite-under-extension}(b)), which contradicts the proper infiniteness of $J^{**}$.

\smnoind
(b) Assume that there is a non-zero finite closed ideal $I\subseteq A$.
By \cite[Corollary 3.6]{CP79}, there is a faithful tracial state $\tau$ on $I$, and it extends to $\tau_0\in T_s(A)$ (see Lemma \ref{lem:finite-under-extension}(b)).
For any  $x\in (I \cap \ker \pi_{\tau_0})_+$, we have $\tau_0(x) = 0$, and the faithfulness of $\tau$ gives $x=0$.
Hence,  $I \cap \ker \pi_{\tau_0} = \{0\}$ and $\ker \pi_{\tau_0}$ is not an essential ideal.
The converse follows directly from part (a).
\end{prf}

\medskip

The above and Remark \ref{rem:main1}(a) give the following.
The reader may note the similarity and difference between Statement (3) below and Proposition \ref{prop:pure-inf>III}(a).

\medskip

\begin{cor}\label{cor:cp-type-III-anti-fin}
The following statements are equivalent.
\begin{enumerate}
\item $A$ is of type $\ttr$.
\item All hereditary $C^*$-subalgebras of $A$ are essentially infinite.
\item For every hereditary $C^*$-subalgebra $B\subseteq A$, there is an essential ideal $I$ of $B$ with $I^{**}$ being a properly infinite $W^*$-algebra.
\end{enumerate}
In particular, all type $\ttr$ $C^*$-algebra are essentially infinite.
\end{cor}

\medskip

\begin{eg}\label{eg:essent-inf}
(a) If $G$ is a separable connected non-amenable locally compact group, then it follows from Proposition \ref{prop:equiv-anti-fin}(a) and the main result of \cite{Ng-str-amen-rep} that the reduced group $C^*$-algebra $C^*_r(G)$ is essentially infinite.

\smnoind
(b) Let $\Gamma$ be a countably infinite amenable group and $c_0(\Gamma)^1$ be the unitalization of $c_0(\Gamma)$.
Consider $\beta$ to be the unique extension of the left translation action of $\Gamma$ on $c_0(\Gamma)$ to $c_0(\Gamma)^1$, and $B:=c_0(\Gamma)^1\rtimes_\beta \Gamma$.
The spectrum $\pr(c_0(\Gamma)^1)$ equals the one-point compactification, $\Gamma \cup \{\infty\}$, of $\Gamma$.
For any $t\in \Gamma\setminus \{e\}$, one has
$$\pr(c_0(\Gamma)^1)^t := \big\{x\in \Gamma \cup \{\infty\}: \ti\beta_t(x) = x \big\} = \{\infty\},$$
where $\ti \beta$ is the induced action of $\Gamma$ on $\pr(c_0(\Gamma)^1)$.
Hence, the action $\beta$ is almost free in the sense of \cite[Definition 4]{LN}.
Suppose that $J$ is a non-zero closed ideal of $B$.
It follows from \cite[Theorem 9]{LN} that there is a non-zero $\beta$-invariant ideal $I\subseteq c_0(\Gamma)^1$ with $I\rtimes_\beta \Gamma \subseteq J$.
Since $c_0(\Gamma)$ is essential in $c_0(\Gamma)^1$, we know that $I\cap c_0(\Gamma)\neq \{0\}$, and this gives
$$I\rtimes_\beta \Gamma \cap c_0(\Gamma)\rtimes_\beta \Gamma \neq \{0\}.$$
Consequently, $\CK(\ell^2(\Gamma)) = c_0(\Gamma)\rtimes_\beta \Gamma$ is an essential ideal of $B$, and Proposition \ref{prop:equiv-anti-fin}(a) tells us that $B$ is essentially infinite.

On the other hand, the short exact sequence $0\to c_0(\Gamma) \to c_0(\Gamma)^1 \to \BC\to 0$ induces the exact sequence
$$0\to c_0(\Gamma)\rtimes_\beta \Gamma \to B \to C^*(\Gamma)\to 0.$$
From which, we know that there exists a tracial state on $B$.

\smnoind
(c) Those simple $C^*$-algebras obtained by R\o rdam in \cite{Rord03} are essentially infinite.
In fact, a simple $C^*$-algebra can only be either finite or essentially infinite (see Remark \ref{rem:fin-ess-inf}(d)).
If $A$ is one of the $C^*$-algebras in \cite{Rord03}, as $A$ contains an infinite projection $p$, this algebra cannot be finite (otherwise, for any proper subprojection $q\leq p$ one obtains, by  \cite[Theroem 3.4]{CP79}, $\tau\in T_s(A)$ with $\tau(p-q) > 0$, which prevents $p\sim q$).
\end{eg}

\medskip

It was said in \cite{Rord03} that the non-type $\mathrm{I}$ algebras obtained there should be $C^*$-analogues of non-type $\tw_1$ and non-type $\ttr$ factors.
In the framework here, those algebras cannot be discrete (since simple discrete $C^*$-algebras are of the form $\CK(H)$; see \cite[Corollary 4.5]{NW-Mv-class} or \cite[Corollary 2.4]{PZ00}).
However, we do not know whether they are type $\tw$ essentially infinite or type $\ttr$.
Notice that a projection in a $C^*$-algebra being finite in the sense of Murray-von Neumann is much weaker than it being finite in the sense of Cuntz-Pedersen (e.g., the identity of any unital $C^*$-algebra without a non-trivial projection is finite in the sense of Murray-von Neumann).

%\medskip

%We end this section with the following characterization of type $\ttr$ that is similar to Proposition \ref{prop:equiv-anti-fin}(a).
%As in \cite{CP79}, we denote by $\T(A)$ the set of lower semi-continuous semi-finite trace on $A$.
%The following lemma follows from the argument of \cite[Lemma 4.6]{CP79}.

%\medskip

%\begin{lem}
%If $J$ is any ideal of $A$, then for any $\varphi\in \T(J)$,
%$\varphi_0(a) := \sup \{\varphi(y): y\in J, \varphi(y) < +\infty \text{ and } y\leq b  \text{ for some } b\sim a\}$
%defines an element in $\T(A)$. Q: is $\varphi_0$ semi-finite?
%\end{lem}

%\medskip

%\begin{prop}\label{prop:equiv-ttr}
%$A$ is type $\ttr$ if and only if there is an essential ideal $I$ of $A$ such that for any non-zero hereditary $C^*$-subalgebra $B\subseteq I$, one has $T_s(B) = \emptyset$.
%\end{prop}
%\begin{prf}
%$\Rightarrow)$.
%For each $\varphi\in \T(A)$, the set $L_\varphi:= \{a\in A: \varphi(a^*a) = 0\}$ is a closed ideal of $A$.
%We put $I:= \bigcap_{\varphi\in \T(A)} L_\varphi$.
%Suppose on the contrary that $I^\bot \neq (0)$ and take any $x\in I^\bot_+\setminus \{0\}$.
%If $y\in I^\bot_+$ such that $y\leq x$ and $y\sim x$, then

%\end{prf}

\medskip

\section{The classification scheme in terms of discreteness, type $\tw$ and type $\ttr$}

\medskip

In this section, we will consider the remaining statements in Theorem \ref{thm:main3}.
Let us begin with the following two simple lemmas.

\medskip

\begin{lem}\label{lem:ideal-non-discrete-alg}
If $A$ is not discrete, then $A$ contains  a
non-zero closed ideal of  either type $\tw$ or type $\ttr$.
\end{lem}
\begin{prf}
By Theorem \ref{thm:sf-C-st-alg}(a), there is a non-zero anti-liminary closed ideal $I\subseteq A$.
If $I$ does not contain a non-zero finite element as well, then it is of type $\ttr$.
Otherwise, $I$ will contain a non-zero finite hereditary $C^*$-subalgebra $B$ (see Lemma \ref{lem:finite}(c)).
The non-zero ideal $J := \overline{IBI}$ of $I$ is semi-finite, because of Remark \ref{rem:str-Mori-equiv}(d)
and Theorem \ref{thm:sf-C-st-alg}(c).
Hence, $J$ is of type $\tw$.
\end{prf}

\medskip

\begin{lem}\label{lem:sum-ideal-essent}
Let $I$ and $J$ be two closed ideals of $A$ such that $I\cap J = (0)$ and $I+J$ is essential.
If $\varphi:A\to M(I)$ is the canonical $^*$-homomorphism, then $\varphi$ restricts to an injection on $J^\bot$.
\end{lem}
\begin{prf}
If $J^\bot \cap I^\bot\neq (0)$, then $J^\bot \cap I^\bot \cap (I+J) \neq (0)$,
	and Proposition \ref{prop:fact-open-hered}(c) gives the contradiction that either $J^\bot \cap I^\bot \cap I\neq (0)$ or $J^\bot \cap I^\bot \cap J \neq (0)$.
	Consequently, $\varphi|_{J^\bot}$ is injective (as $\ker \varphi = I^\bot$).
\end{prf}

\medskip

\begin{thm}\label{thm:decomp}
Let $A$ be a $C^*$-algebra.

\smnoind
(a) There exists the largest discrete (respectively, semi-finite, type $\tw$ and type $\ttr$) hereditary $C^*$-subalgebra $A_\dc$ (respectively, $A_\SF$, $A_\tw$ and $A_\ttr$) of $A$.
Moreover, $A_\dc$, $A_\SF$, $A_\tw$ and $A_\ttr$ are ideals of $A$ such that $A_\dc$, $A_\tw$ and $A_\ttr$ are disjoint and $A_\ttr\cap A_\SF = (0)$.

\smnoind
(b) $A_\dc + A_\tw + A_\ttr$ is an essential ideal of $A$ and $A_\dc + A_\tw$ is an essential ideal of $A_\SF$.

\smnoind
(c) $A_\al:=A_\dc^\bot$ is the largest anti-liminary hereditary $C^*$-subalgebra of $A$.

\smnoind
(d) $A_\dc = A_\al^\bot = (A_\tw + A_\ttr)^\bot$, $A_\tw = A_\SF\cap A_\dc^\bot$, $A_\dc = A_\SF\cap A_\tw^\bot$, $A_\ttr = A_\SF^\bot =  (A_\dc + A_\tw)^\bot$ and $A_\SF=A_\ttr^\bot$.

\smnoind
(e) If $A$ is a $W^*$-algebra, then $A_\dc$, $A_\tw$ and $A_\ttr$ are respectively, the type {\tI}, the type $\tw$ and the type $\ttr$ $W^*$-algebra
summands of $A$.
\end{thm}
\begin{prf}
(a) Let us first construct the largest type $\tw$ hereditary $C^*$-subalgebra of $A$ and verify that it is an ideal.
In order to do so, we will show that the collection $\CJ_\tw$ of all type $\tw$ closed ideals of $A$ is a directed set, and that 
$$A_\tw\ :=\ \overline{{\sum}_{J\in
		\CJ_\tw} J},$$
is the largest type $\tw$ hereditary $C^*$-subalgebra of $A$, which is clearly an ideal.

For the first claim, we will verify that $J_1+J_2\in \CJ_\tw$ for any arbitrary elements $J_1, J_2 \in \CJ_\tw$.
In fact, if $B$ is a non-zero
hereditary $C^*$-subalgebra of $J_1+J_2$, then by Proposition
\ref{prop:fact-open-hered}(c), we have $B\cap J_1\neq (0)$ or $B\cap J_2\neq (0)$.
If $B\cap J_1\neq (0)$, then one obtains a non-zero finite element $x\in \RF^{B\cap J_1}\subseteq \RF^B$.
Similarly, if $B\cap J_2\neq (0)$, then $B$ will also contains a non-zero finite element.
On the other hand, suppose that there exists a non-zero
abelian element $a \in (J_1+J_2)_+$.
Set $D := \overline{a(J_1+J_2)a}$.
Again, it follows from Proposition
\ref{prop:fact-open-hered}(c) that $D\cap J_1\neq (0)$ or $D\cap J_2\neq (0)$.
However, this contradicts with $J_1, J_2\in \CJ_\tw$, because $D\cap J_1\subseteq J_1$  and $D\cap J_2\subseteq J_2$ are abelian hereditary $C^*$-subalgebras.
Thus, $J_1 + J_2\in \CJ_\tw$.

Now, for any nonzero closed ideal $J$ of $A$, let us denote by $p_J$ the element in  $\op(A)\cap Z(A^{**})$ with $J = \her_A(p_J)$.
It is easy to see that $p_{A_\tw}$ is the $\sigma(A^{**}, A^*)$-limit of the increasing net $\{p_J\}_{J\in \CJ_\tw}$.
Suppose that $e\in \op(A)$  satisfying $\her_A(e)\subseteq
A_\tw$.
Then $$
e\ =\ e p_{A_\tw} e\ =\ w^*\text{-}{\lim}_{J\in \CJ_\tw} e p_Je,
$$
and there is $J\in \CJ_\tw$ with $ep_J \neq 0$.
This means that
$\her_A(e)\cap J$ is non-zero (see Proposition \ref{prop:fact-open-hered}(b)), and  is a hereditary $C^*$-subalgebra of $J$. 
Hence, $\her_A(e)\cap J$ contains a
non-zero finite element.
On the other hand, assume that there is a non-zero abelian positive element $a$ in $A_\tw$.
Consider $e\in \op(A)$ such that $\her_A(e) = \overline{aAa} \subseteq A_\tw$.
As in the above, one can find $J\in \CJ_\tw$ with $ep_J \neq 0$.
Thus, $\her_A(e)\cap J$ is a non-zero abelian hereditary $C^*$-subalgebra of $J$, which contradicts $J$ being of type $\tw$.
Consequently, one has $A_\tw\in \CJ_\tw$.
Finally, if $B\subseteq A$
is a hereditary $C^*$-subalgebra of type $\tw$, then, by Remark \ref{rem:str-Mori-equiv}(d) and Theorem \ref{thm:sf-C-st-alg}(c), one knows that $B\subseteq \overline{ABA}\subseteq A_\tw$.

The existences of $A_\SF$, $A_\tw$ and $A_\ttr$ as well as the fact that they are ideals follow from similar arguments.
The disjointness statements are obvious.

\smnoind (b)
This part follows directly from Lemma
\ref{lem:ideal-non-discrete-alg} (namely, every non-zero non-discrete closed ideal interests  $A_\tw$ or $A_\ttr$).

\smnoind
(c) If $A_\al$ contains a non-zero abelian positive element $x$, then $\overline{xAx}\subseteq A_\dc$, and we have a contradiction that $x\in A_\al\cap A_\dc$.
Thus, $A_\al$ is anti-liminary.
Furthermore, if $B\subseteq A$ is a non-zero anti-liminary hereditary $C^*$-subalgebra, then $BA_\dc B = A_\dc \cap B = (0)$, which means that $B\subseteq A_d^\bot = A_\al$.

\smnoind
(d) By parts (a) and (c), one has $A_\dc\subseteq A_\al^\bot \subseteq (A_\tw+A_\ttr)^\bot$.
Since $A_\dc\cap (A_\tw+A_\ttr) = (0)$ and $A_\dc + A_\tw + A_\ttr$ is
an essential ideal of $A$, we obtain from Lemma \ref{lem:sum-ideal-essent} a $^*$-monomorphism from
$(A_\tw+A_\ttr)^\bot$ to $M(A_\dc)$ whose image clearly contains $A_\dc$.
Therefore, Remark \ref{rem:essetial}(d) and Proposition \ref{prop:type-hered}(b) tells us that $(A_\tw+A_\ttr)^\bot$ is discrete; i.e, $(A_\tw+A_\ttr)^\bot\subseteq A_\dc$.
These give the first equality. 

Secondly, since $A_\tw A_\dc  = (0)$ and $A_\dc + A_\tw$ is essential in $A_\SF$ (see part (b)), we have $A_\tw \subseteq A_\SF\cap A_\dc^\bot$ and Lemma \ref{lem:sum-ideal-essent} produces a $^*$-monomorphism from $A_\SF\cap A_\dc^\bot$ to $M(A_\tw)$.
As in the above, this gives the opposite inclusion $A_\SF\cap A_\dc^\bot\subseteq A_\tw$.
The equality $A_\dc = A_\SF\cap A_\tw^\bot$ follows from the same argument.

To establish the fourth equality, let us set $J_3:= (A_\dc + A_\tw)^\bot$.
As $A_\ttr A_\SF  = (0)$ and $A_\dc + A_\tw\subseteq A_\SF$,  we see that $A_\ttr \subseteq A_\SF^\bot \subseteq J_3$.
It remains to show that $J_3$ is of type $\ttr$.
In fact, suppose on the contrary that $J_3$ contains a non-zero finite hereditary $C^*$-subalgebra $B$. 
Then $B\cap A_d = (0)$ (because $J_3\subseteq A_\dc^\bot$).
This implies that $B$ is of type $\tw$, which contradicts $J_3\subseteq A_\tw^\bot$.

Finally, it is clear that $A_\SF\subseteq A_\ttr^\bot$.
Conversely, since $A_\SF \cap A_\ttr = (0)$ and $A_\SF + A_\ttr$ contains the essential ideal $A_I + A_\tw + A_\ttr$ of $A$, one learns from Lemma \ref{lem:sum-ideal-essent} that there is a $^*$-monomorphism from
$A_\ttr^\bot$ to $M(A_\SF)$ whose image clearly contains $A_\SF$.
It now follows from Remark \ref{rem:essetial}(d) and Proposition \ref{prop:type-hered}(b) that $A_\ttr^\bot$ is semi-finite, and thus $A_\ttr^\bot\subseteq A_\SF$.

\smnoind
(e) This part follows from Proposition \ref{prop:RR0}(d) and the fact that $A_\dc$, $A_\tw$ and $A_\ttr$ are $\sigma(A,A_*)$-closed (because of part (d)).
\end{prf}

\medskip

This theorem produces the following result concerning the stability of the above decomposition under strong Morita equivalence (see Remark \ref{rem:str-Mori-equiv}(c) for the notation).

\medskip

\begin{cor}\label{cor:morita-equiv-decomp}
If $X$ is a full Hilbert $A$-module and $B:= \CK_A(X)$, then
$B_\dc = \CK_{A}(XA_\dc)$, $B_\tw= \CK_{A}(XA_\tw)$ and $B_\ttr= \CK_{A}(XA_\ttr)$.
\end{cor}
\begin{prf}
It follows from Remark \ref{rem:str-Mori-equiv}(c) that the assignment $I\mapsto \CK_A(XI)$ is a bijection from the set of closed ideals of $A$ to that of $B$, and that the ideal $I$ of $A$ is strongly Morita equivalent to $\CK_A(XI)$.
Since $A_\dc$ (respectively, $A_\tw$ and $A_\ttr$) is the largest discrete (respectively, type $\tw$ and type $\ttr$) ideal of $A$, we know from the above and Theorem \ref{thm:sf-C-st-alg}(c) that $B_\dc$ (respectively, $B_\tw$ and $B_\ttr$) is the largest discrete (respectively, type $\tw$ and type $\ttr$) ideal of $B$.
\end{prf}

\medskip

On the other hand, we have the following result. 
The meaning of universal $\#$ normal quotient in its statement is understood in similar ways as those of parts (a) and (b) of Theorem \ref{thm:inf-part}.

\medskip

\begin{cor}\label{cor:quot}
(a) $A/A_\al$ (respectively, $A/A_\dc$, $A/A_\ttr$ and $A/A_\SF$)
is the universal discrete (respectively, anti-liminary, semi-finite and type $\ttr$) normal quotient of $A$.

\smnoind
(b)  $A_\SF/A_\tw$ (respectively, $A_\SF/A_\dc$) is the universal discrete (respectively, type $\tw$) normal quotient of $A_\SF$.
\end{cor}
\begin{prf}
(a) It follows from Theorem \ref{thm:decomp}(c) (respectively, Theorem \ref{thm:decomp}(d)) that $A/A_\al = A/A_\dc^\bot$ (respectively, $A/A_\dc = A/A_\al^\bot$, $A/A_\ttr = A/A_\SF^\bot$ and $A/A_\SF = A/A_\ttr^\bot$), and hence it is a normal quotient.  
Moreover, we know from Proposition \ref{prop:type-hered}(c) that $A/A_\al$ is discrete (respectively, $A/A_\dc$ is anti-liminary, $A/A_\ttr$ is semi-finite and $A/A_\SF$ is of type $\ttr$).

Let $I$ be an ideal of $A$ and $q:A\to A/I^\bot$ be the quotient map.
Suppose that $A/I^\bot$ is discrete (respectively, anti-liminary, semi-finite and type $\ttr$).
Then $q$ restricts to an injection on $I$ with $q(I)$ being an ideal of $A/I^\bot$.
Proposition \ref{prop:type-hered}(a) tells us that $I$ is discrete (respectively, anti-liminary, semi-finite and type $\ttr$), and so $I\subseteq A_\dc$ (respectively, $I\subseteq A_\al$, $I\subseteq A_\SF$ and $I\subseteq A_\ttr$).
Consequently, $A_\al \subseteq I^\bot$ (respectively, $A_\dc \subseteq I^\bot$, $A_\ttr \subseteq I^\bot$ and $A_\tw \subseteq I^\bot$),  because of parts (c) and (d) of Theorem \ref{thm:decomp}, and we conclude that $q$ factors through $A/A_\al$ (respectively, $A/A_\dc$, $A/A_\ttr$ and $A/A_\SF$).

\smnoind
(b) This part follows from the argument of part (a).
Notice that $A_\tw = A_\SF\cap A_\dc^\bot$ and $A_\dc = A_\SF\cap A_\tw^\bot$ (see Theorem \ref{thm:decomp}(d)).
\end{prf}

\medskip

Furthermore, we want to see how type decompositions pass to hereditary $C^*$-subalgebras, essential extensions and normal quotients.

\medskip

\begin{prop}\label{prop:decomp}
Let $B\subseteq A$ be a hereditary $C^*$-subalgebra and $I, J\subseteq A$ be closed ideals with $J$ being essential.
Let $q: A \to A/I^\bot$ be the quotient map.

\smnoind
(a) $B_\dc = A_\dc\cap B$, $B_\SF = A_\SF\cap B$, $B_\al = A_\al\cap B$,  $B_\tw = A_\tw\cap B$ and $B_\ttr = A_\ttr\cap B$.

\smnoind
(b) Consider $\#= \dc, \SF, \al, \tw$ or $\ttr$. 
Then $q(A_\#)$ is an essential ideal of $(A/I^\bot)_\#$ and 
\begin{equation}\label{eqt:A/I-bot}
(A/I^\bot)_\# = q(\{a\in A: aI\subseteq I_\#\}).
\end{equation}
In particular, $A_\# = \{a\in A:aJ \subseteq J_\#\}$.

\smnoind
(c) $A_\tw = \{x\in A: xJ_\dc = (0) \text{ and } xJ \subseteq J_\SF\}$

\smnoind
(d) $A_\ttr = \{x\in A: xJ_\SF = (0)\} = \{x\in A: xJ_\dc = (0)  \text{ and } xJ_\tw = (0)\}$.
\end{prop}
\begin{prf}
(a) By Proposition \ref{prop:type-hered}(a), $A_\dc\cap B$, $A_\SF\cap B$, $A_\al\cap B$,  $A_\tw\cap B$ and $A_\ttr\cap B$ are respectively, discrete, semi-finite, anti-liminary, type $\tw$ and type $\ttr$ hereditary $C^*$-subalgebras of $B$.
On the other hand, if $D\subseteq B$ is a hereditary $C^*$-subalgebra which is discrete (respectively, semi-finite, anti-liminary, type $\tw$ and type $\ttr$), then $D$ is a hereditary $C^*$-subalgebra of $A$ and hence is contained in $A_\dc$ (respectively, $A_\SF$, $A_\al$, $A_\tw$ and $A_\ttr$).
These give the required statement.

\smnoind
(b) Let us first establish the third claim; i.e., $A_\# = \check A_\#$, where
$$\check A_\#\ :=\ \{x\in A: xJ \subseteq J_\#\} \qquad (\# = \dc, \SF, \al, \tw, \ttr).$$
We will consider only the case when $\# = \tw$ (because the other cases follow from similar arguments).
In fact, it follows from  $A_\tw J = J_\tw$ (see part (a)) that $A_\tw \subseteq \check A_\tw$.
Consider $B\subseteq \check A_\tw$ to be an arbitrary non-zero hereditary $C^*$-subalgebra.
Then $B\cap J = BJB$ is contained in $J_\tw$ and $BJB\neq (0)$, as $J$ is essential.
Thus, we obtain a non-zero finite element in $(B\cap J)_+$.
On the other hand, suppose that $\check A_\tw$ contains a non-zero abelian positive  element $a$.
If we set $B := \overline{aA a}$, then $B\cap J$ is a non-zero abelian hereditary $C^*$-subalgebra of $J_\tw$, which is absurd.
Consequently, $\check A_\tw$ is of type $\tw$ and is contained in $A_\tw$.

Secondly, we will verify Relation \eqref{eqt:A/I-bot}.
Note that, by considering the canonical $^*$-monomorphism $\varphi:A/I^\bot\to M(I)$, one sees that $\varphi(A/I^\bot)$ as an essential extension of $I$ (see Remark \ref{rem:essetial}(d)).
Thus, the statement concerning $J_\#$ implies Relation \eqref{eqt:A/I-bot}, since we have
\begin{equation*}
\varphi((A/I^\bot)_\#)
\ = \ \{\varphi(x): x\in A/I^\bot; \varphi(x)I \subseteq I_\#\}
\ = \ \{\varphi(q(a)): a\in A; \varphi(q(a))I \subseteq I_\#\}
\end{equation*}
as well as $\varphi(q(a))x = ax$ ($a\in A; x\in I$). 

Thirdly, we will show that $q(A_\#)$ is an essential ideal of $(A/I^\bot)_\#$.
Indeed, by Relation \eqref{eqt:orth-comp} and part (a) above, one has
$$I^\bot \cap A_\dc = (I_\dc)^\bot \cap A_\dc,$$
which produces a $^*$-monomorphism from $q(A_\dc) \cong A_\dc / (I^\bot \cap A_\dc)$ to $M(I_\dc)$ with its image containing $I_\dc$.
Therefore, Remark \ref{rem:essetial}(d) and Proposition \ref{prop:type-hered}(b) imply that the ideal $q(A_\dc)$ of $A/I^\bot$ lies inside $(A/I^\bot)_\dc$.
Moreover, suppose that $J\subseteq (A/I^\bot)_\dc$ is a non-zero closed ideal.
As $q$ restricts to an injection from $q^{-1}(J)\cap I^{\bot\bot}$ to an ideal of the discrete $C^*$-algebra $(A/I^\bot)_\dc$, we know that $q^{-1}(J)\cap I^{\bot\bot} \subseteq A_\dc$ (because of Proposition \ref{prop:type-hered}(a)) and that $q^{-1}(J)\cap I^{\bot\bot}\neq \{0\}$ (since $q^{-1}(J)\nsubseteq I^\bot$).
These show that $J\cap q(A_\dc)\neq \{0\}$.
The arguments for the statements concerning $q(A_\SF)$, $q(A_\al)$, $q(A_\tw)$ and $q(A_\ttr)$ are the same.

\smnoind
(c) Obviously, for any $x\in A$, one has $xJ_\dc = (0)$ if and only if $xJJ_\dc = (0)$.
Thus, this part follows from part (b) and the fact that $J_\tw = J_\SF \cap J_d^\bot$ (see Theorem \ref{thm:decomp}(d)).

\smnoind
(d) Note that we also have $xJ_\# = (0)$ if and only if $xJJ_\# = (0)$ for $\# = \tw, \SF$.
Hence, this part follows from part (b) as well as the two equalities $J_\ttr = J_\SF^\bot = (J_\dc + J_\tw)^\bot$
(as given in Theorem \ref{thm:decomp}(d)).
\end{prf}

\medskip

\begin{cor}\label{cor:decomp}
Let $A_{\dc,1} := A_\dc\cap A_1$, $A_{\dc,\infty} := A_\dc\cap A_\infty$,
$A_{\tw,1} := A_\tw\cap A_1$ and $A_{\tw,\infty} := A_\tw\cap A_\infty$.
Then $A_{\dc,1}$, $A_{\dc,\infty}$, $A_{\tw,1}$ and $A_{\tw,\infty}$ are the largest discrete finite ideal, the largest discrete essentially infinite ideal, the largest type $\tw$ finite ideal and the largest type $\tw$ essentially infinite ideal respectively.
These are disjoint normal ideals.
Moreover, $A_{\dc,1} + A_{\dc,\infty} + A_{\tw,1} + A_{\tw,\infty} + A_\ttr$ is an essential ideal of $A$, or equivalently,
$$A_{\dc,1} \oplus A_{\dc,\infty} \oplus A_{\tw,1} \oplus A_{\tw,\infty} \oplus A_\ttr \subseteq A \subseteq M(A_{\dc,1}) \oplus M(A_{\dc,\infty}) \oplus M(A_{\tw,1}) \oplus M(A_{\tw,\infty}) \oplus M(A_\ttr).$$
\end{cor}
\begin{prf}
We first note that by Theorem \ref{thm:inf-part}(d),
\begin{equation}\label{eqt:AIIinf}
A_{\dc,1} = (A_\dc)_1, \quad
A_{\dc, \infty} = (A_\dc)_\infty, \quad A_{\tw, 1} = (A_\tw)_1 \quad \text{and} \quad A_{\tw, \infty} = (A_\tw)_\infty.
\end{equation}
Thus, the first claim follows from Proposition \ref{prop:type-hered}(a), parts (a)  and (b) of Theorem \ref{thm:inf-part} as well as Theorem \ref{thm:decomp}(a).
These four ideals are clearly disjoint and they are normal because of Remark \ref{rem:fin-ess-inf}(a), Theorem \ref{thm:inf-part}(a), parts (c) and (d) of Theorem \ref{thm:decomp} as well as the fact that $I^\bot \cap J^\bot = (I+J)^\bot$ for any two closed ideals $I,J\subseteq A$.
In order to establish the third statement, we consider a non-zero closed ideal $J\subseteq A$.
By Theorem \ref{thm:decomp}(b) and Proposition \ref{prop:fact-open-hered}(c), we know that
$$J\cap A_\dc\neq \{0\}, \quad J\cap A_\tw\neq \{0\} \quad \text{or} \quad J\cap A_\ttr\neq \{0\}.$$
If $J\cap A_\dc\neq \{0\}$, we know from Remark \ref{rem:fin-ess-inf}(b) as well as  the first two equalities of \eqref{eqt:AIIinf} that $J\cap A_{\dc,1}\neq \{0\}$ or $J\cap A_{\dc,\infty}\neq \{0\}$.
Similarly,  if $J\cap A_\tw\neq \{0\}$, then $J\cap A_{\tw,1}\neq \{0\}$ or $J\cap A_{\tw,\infty}\neq \{0\}$.
Thus, $A_{\dc,1} + A_{\dc,\infty} + A_{\tw,1} + A_{\tw,\infty} + A_\ttr$ is an essential ideal.
The displayed relations follow from the well-known fact that $M(B\oplus D) = M(B) \oplus M(D)$ for two $C^*$-algebras $B$ and $D$.
\end{prf}

\medskip

Observe that the $C^*$-algebras $M(A_{\dc,1})$, $M(A_{\dc,\infty})$, $M(A_{\tw,1})$, $M(A_{\tw,\infty})$ and $M(A_\ttr)$ are also discrete finite, discrete essentially infinite , type $\tw$ finite, type $\tw$ essentially infinite and type $\ttr$, respectively (by Proposition \ref{prop:type-hered}(b) and the argument of Theorem \ref{thm:inf-part}).
Furthermore, one may obtain results similar to those in Corollary \ref{cor:quot} for $A_{\dc,1}$, $A_{\dc,\infty}$, $A_{\tw,1}$ and $A_{\tw,\infty}$, but we leave them to the readers.

\bigskip

\section{Some special cases}

In this section, we will consider special classes of $C^*$-algebras for which we have a better decomposition in Corollary \ref{cor:decomp}.

\medskip

\subsection{Prime $C^*$-algebras}

The first special class is that of prime $C^*$-algebras.
Recall that a $C^*$-algebra $A$ is \emph{prime} if $\{0\}$ is a prime ideal of $A$, or equivalently, $A$ has no non-trivial normal ideal.
Since the ideals $A_{\dc,1}$, $A_{\dc,\infty}$, $A_{\tw,1}$, $A_{\tw,\infty}$ and $A_\ttr$ are normal (see Corollary \ref{cor:decomp}), we have the following.

\medskip

\begin{prop}\label{prop:prime>5types}
Any prime $C^*$-algebra is of one of the five types: discrete finite, discrete essentially infinite, type $\tw$ finite, type $\tw$  essentially infinite or type $\ttr$.	
\end{prop}

\medskip

It is well-known that the kernel of a factor representation of a $C^*$-algebra is a prime ideal.
If the $C^*$-algebra is separable, one has the strong converse that every prime ideal is primitive.
These give the following well-known fact for a $C^*$-algebra $A$:

\begin{enumerate}[(P1).]
\item If there is an injective $^*$-homomorphism $\varphi$ from $A$ to a factor $M$ with $\varphi(A)'' = M$, then $A$ is prime.

\item If $A$ is a separable prime $C^*$-algebra, there exists a faithful irreducible $^*$-representation of $A$.
\end{enumerate}

\medskip

Clearly, a $W^*$-algebra is a factor if and only if it is a prime $C^*$-algebra.
Because of this, as well as Proposition \ref{prop:quasi-Ston}(b) and Corollary \ref{cor:ess-simp} below, one may regard prime $C^*$-algebras as ``$C^*$-algebra factors''.

\medskip

Let us present the following description of discrete prime $C^*$-algebras.
This result is more or less a reformation of \cite[Corollary 2.4]{PZ00}.

\medskip

\begin{prop}\label{prop:prime-discrete}
Let $A$ be a $C^*$-algebra.

\smnoind
(a) The following statements are equivalent.
\begin{enumerate}
	\item $A$ is prime and discrete.
	\item There exists a Hilbert space $\KH$ such that $\CK(\KH)\subseteq A \subseteq \CL(\KH)$.
	\item $A$ contains $\CK(\KH)$ as an essential ideal for some Hilbert space $\KH$.
\end{enumerate}

\smnoind
(b) $A$ is prime, discrete and finite if and only if $A\cong M_n(\BC)$ for some $n\in \BN$.
\end{prop}
\begin{prf}
(a) $(1)\Rightarrow (2)$.
This follows directly from \cite[Corollary 2.4]{PZ00}.

\smnoind
$(2)\Rightarrow (3)$.
This follows from Remark \ref{rem:essetial}(d).

\smnoind
$(3)\Rightarrow (1)$.
Since $A$ contains $\CK(\KH)$ as an essential ideal, Proposition \ref{prop:type-hered}(b) implies that $A$ is discrete.
Moreover, if $J\subseteq A$ is a non-zero closed ideal, then $J \cap \CK(\KH) = \CK(\KH)$ (as $\CK(\KH)$ is simple) which means that $J^\bot \subseteq \CK(\KH)^\bot = \{0\}$.
In other words, $A$ has no non-zero normal ideal and hence it is prime.

\smnoind
(b) It is well-known that $M_n$ is a discrete finite simple $C^*$-algebra.
Conversely, if $\KH$ is as in Statement (2) of part (a), then the ideal $\CK(\KH)$ of $A$ is finite.
Hence, $n:= \dim \KH < \infty$ and $A\cong  M_n(\BC)$.
\end{prf}

\medskip

It follows that the unitalization of $\CK(\ell^2)$ is a discrete essentially infinite prime $C^*$-algebra.
In the following, we give some more examples of prime $C^*$-algebras of different types which are not $W^*$-algebra.

\medskip

\begin{eg}\label{eg:prime}
(a) If $\Gamma$ is a countably infinite amenable group, and $\beta$ is as in Example \ref{eg:essent-inf}(b), then $c_0(\Gamma)^1\rtimes_\beta \Gamma$ contains $\CK(\ell^2(\Gamma))$ as an essential ideal, and hence it is a discrete  essentially infinite prime $C^*$-algebra.

\smnoind
(b) Let $\Gamma$ be any countable ICC group.
Since the group von Neumann algebra $L(\Gamma)$ is a factor, Statement (P1)  above implies that the reduced group $C^*$-algebra $C^*_r(\Gamma)$ is prime.
Moreover, since $C^*_r(\Gamma)$ has a faithful trace, it is finite.
Consequently, by Propositions \ref{prop:prime>5types} and Proposition \ref{prop:prime-discrete}(b), $C^*_r(\Gamma)$ is type $\tw$ finite.

As in the case of $W^*$-algebra, if, in addition, $\Gamma$ is amenable, then $C^*_r(\Gamma)$ is nuclear, and is non-isomorphic to the non-nuclear type $\tw$ finite prime $C^*$-algebra $C^*_r(\mathbb{F}_2)$, where $\mathbb{F}_2$ is the free group on two generators. 

\smnoind
(c) Suppose that $A\ncong \CK(\ell^2)$ is a simple $AF$-algebra that does not have a tracial state.
Then $A$ is a type $\tw$ essentially infinite prime $C^*$-algebra.
In fact, $A$ is not finite as it has no tracial state.
By \cite[Proposition 4.11]{CP79}, $A$ is semi-finite.
However, since discrete simple $C^*$-algebras are of the form $\CK(\ell^2)$, we know that $A$ is not discrete.

\smnoind
(d) The Calkin algebra $\CL(\ell^2)/\CK(\ell^2)$ is a simple purely infinite $C^*$-algebra.
Hence it is a type $\ttr$ prime $C^*$-algebra (by Proposition \ref{prop:pure-inf>III}(b)).
\end{eg}

\medskip

\subsection{$C^*$-algebras with extremely disconnected primitive spectrum}
The second special class are $C^*$-algebras of which all normal ideals are complementary.
Before looking at this class, let us first give the following result concerning normal ideals of essential extensions.

\medskip

\begin{lem}\label{lem:ann-id-ess-ext}
If $A$ is an essential extension of a $C^*$-algebra $B$, the map $\Gamma:J\mapsto J\cap B$ is a bijection from the set of normal ideals of $A$ onto that of $B$.
\end{lem}
\begin{prf}
Let $I\subseteq A$ be a closed ideal.
We first show that
\begin{equation}\label{eqt:cap-bot}
I^\bot = (I\cap B)^\bot.
\end{equation}
Indeed, one clearly has
$I^\bot \subseteq (I\cap B)^\bot$.
Suppose on the contrary that there is $x\in (I\cap B)^\bot \setminus I^\bot$.
Then $xI\neq (0)$ and hence $\overline{IxI}$ is a non-zero ideal of $A$.
As $B$ is an essential ideal of $A$, one can find $b\in \overline{IxI}\cap B\setminus \{0\}$.
Thus, $b$ can be approximated in norm by elements of the form $\sum_{k=1}^n y_kxz_k$ where $y_1,...,y_n, z_1,...,z_n\in I$, and $bb^*$ can be approximated by elements of the form $\sum_{k=1}^n y_kxz_kb^*$.
Since $z_kb^*\in I\cap B$, we obtain the contradiction that $bb^* = 0$.

Now, if $J = I^\bot$, then by Relation \eqref{eqt:orth-comp}, one has $J\cap B = (I\cap B)^\bot \cap B$ and hence $\Gamma(J)$ is a normal ideal of $B$.
Furthermore, as
\begin{equation}\label{eqt:J-cap-A-cl}
\Gamma(J)^{\bot\bot} = (J\cap B)^{\bot \bot} = J^{\bot\bot} = J
\end{equation}
(by Relation \eqref{eqt:cap-bot}), we see that $\Gamma$ is injective.

Finally, let $J_0$ be a normal ideal of $B$ and  $I_0\subseteq B$ be an ideal satisfying $J_0= I_0^\bot \cap B$.
If $J:= I_0^\bot$, then obviously, $\Gamma(J) = J_0$.
\end{prf}

\medskip

Recall that a topological space (not necessarily Hausdorff) is \emph{extremely disconnected} if the closure of every open subset is again open.
It happens that those $C^*$-algebras described in the beginning of this subsection are precisely those with extremely disconnected primitive spectra
(equipped with the hull-kernel topology).
We will   give some other equivalent forms of these algebras in the following proposition (although not all of them are needed in this paper).

\medskip

\begin{prop}\label{prop:cent-AW-st}
The following statements are equivalent.
\begin{enumerate}
\item The primitive spectrum $\pr(A)$ is extremely disconnected.
\item $A = I^\bot + I^{\bot\bot}$, for every closed ideal $I\subseteq A$.
\item $\overline{p}^1\in \op(A)$ for every $p\in \op(A)\cap Z(A^{**})$.
\item $\Lambda(q):= qA$ is a surjection (and hence a bijection) from  $\proj(ZM(A))\cup \{0\}$ onto the set of normal ideals of $A$.
\item $\pr(M(A))$ is extremely disconnected.
\item $ZM(A)$ is an $AW^*$-algebra and $\Delta: J\mapsto J\cap ZM(A)$ is an injection (equivalently, a bijection) from the set of normal ideals of $M(A)$
to the set of normal ideals  of $ZM(A)$.
\item $A$ is boundedly centrally closed (in the sense of \cite[p.165]{AM94}).
\end{enumerate}
\end{prop}
\begin{prf}
$1)\Leftrightarrow 2)$.
Let $I \subseteq A$ be an ideal.
Since any element $P\in \pr(A)$ is prime, we know that either $I\subseteq P$ or $I^\bot \subseteq P$.
This ensures that
$$I':=\bigcap \pr (A) \setminus \hull (I) \supseteq \bigcap \{P\in \pr(A): I^\bot \subseteq P \} = I^\bot.$$
Conversely, one has $I'\subseteq I^\bot$. 
In fact, suppose on the contrary that $I'I \neq \{0\}$.
Then one can find $P_0\in \pr(A)$ with $I'I \nsubseteq P_0$.
Thus, $P_0\notin \hull(I)$, and we arrive at the contradiction that $I'\subseteq P_0$.
In other words,
$$\overline{\pr (A) \setminus \hull (I)} = \hull (I^\bot).$$

Consequently, $\overline{\pr (A) \setminus \hull (I)}$ is open if and only if $\hull(I^\bot)$ is open, which in turn is equivalent to
$$\pr(A) \setminus \hull(I^\bot) = \overline{\pr(A) \setminus \hull(I^\bot)} = \hull(I^{\bot\bot}).$$
Therefore, $\overline{\pr (A) \setminus \hull (I)}$ being open is the same as
\begin{equation}\label{eqt:hull-J-disj-hull-J-bot}
\hull(I^{\bot})\cap \hull(I^{\bot\bot}) = \emptyset
\end{equation}
(recall that we always have $\overline{\pr (A) \setminus \hull (I^\bot)} = \hull (I^{\bot\bot})$).

Finally, Relation \eqref{eqt:hull-J-disj-hull-J-bot} and $A = I^\bot + I^{\bot\bot}$ are equivalent, because a closed ideal $I_0\subseteq A$ is proper if and only if there is $P\in \pr(A)$ satisfying $I_0\subseteq P$.

\smallskip\noindent
$2)\Rightarrow 3)$.
Suppose that $p\in \op(A)\cap Z(A^{**})$.
As $\her_A(p)^\bot = \her_A(1-\overline{p}^1)$ (see Equality \eqref{eqt:her(e)-bot}), we know from Statement (2) that
$$
\her_A(1-\overline{p}^1) + \her_A(1-\overline{p}^1)^\bot = A.
$$
Hence, for any $a\in A$, there exist $x\in \her_A(1-\overline{p}^1)$ and $y\in \her_A(1-\overline{p}^1)^\bot$ with $a = x + y$.
This shows that $(1-\overline{p}^1)a = x\in A$.
Consequently, $1-\overline{p}^1$ belongs to $M(A)$, and hence is a closed projection of $A$.

\smallskip\noindent
$3)\Rightarrow 4)$.
We first note that $\Lambda$ is always injective, because $\proj(ZM(A))\subseteq \op(A)\cap Z(A^{**})$ and we have Statement (O1).
Moreover, it is clear that $\Lambda(q)$ equals $((1-q)A)^\bot$ and hence is a normal ideal of $A$.
Suppose that $I\subseteq A$ is a closed ideal and $J := I^\bot$.
If $p, q\in \op(A)\cap Z(A^{**})$ with $I = \her_A(p)$ and $J = \her_A(q)$, then Relation \eqref{eqt:her(e)-bot} tells us that $q = 1-\overline{p}^1$.
As $\overline{p}^1$ is both open and closed (by Statement (3)), it belongs to $M(A)$ and so is $q$.
Thus, $J = qA$ as required.

\smallskip\noindent
$4)\Rightarrow 2)$.
If $\Lambda$ is surjective, then for any closed ideal $I\subseteq A$, there is $q\in \proj(ZM(A))\cup \{0\}$ with $I^\bot = qA$, which gives $A = qA + (1-q)A = I^{\bot} + I^{\bot\bot}$.

\smallskip\noindent
$4)\Leftrightarrow 5)$.
By applying the equivalence of Statements (1) and (4) (which was established above) to the $C^*$-algebra $M(A)$, it suffices to show that Statement (4) (for the $C^*$-algebra $A$) is equivalent to
the following statement:
\begin{enumerate}[(1')]
\setcounter{enumi}{3}
\item For any normal ideal $J\subseteq M(A)$, there is $q\in \proj(ZM(A))\cup \{0\}$ with $J = qM(A)$.
\end{enumerate}

In fact, suppose that $J$ is a normal ideal of $M(A)$.
By Lemma \ref{lem:ann-id-ess-ext}, the ideal $J\cap A\subseteq A$ is normal.
Hence, Statement (4) (for the $C^*$-algebra $A$) gives $q\in \proj(ZM(A))\cup\{0\}$ with $J\cap A = qA$.
Observe that if $x\in M(A)$ satisfying $xqA = (0)$, then $x =  (1-q)x$.
This, together with Relation \eqref{eqt:J-cap-A-cl},  tells us that
$$J = (J\cap A)^{\bot\bot} = (qA)^{\bot\bot} = \big((1-q)M(A)\big)^\bot = qM(A)$$
(note that these $\bot$ are taken in $M(A)$).

Conversely, if $I$ is a normal ideal of $A$, then Lemma \ref{lem:ann-id-ess-ext} and Statement (4') produce $q\in \proj(ZM(A))\cup \{0\}$ satisfying $I = qM(A) \cap A = qA$, as is required.

\smallskip\noindent
$1)\Rightarrow 6)$.
Suppose that $\CI$  is a closed ideal of $ZM(A)\cong C_b(\pr(A))$ (by the Dauns-Hofmann theorem; see e.g. \cite[Corollary 4.4.8]{Ped79}).
Then 
$$U_\CI:= \{J\in \pr(A): f(J)\neq 0, \text{ for some }f\in \CI\}$$ 
is open,  and its closure $\overline{U_\CI}$ is an open subset of $\pr(A)$ (by Statement (1)).
One may regard $\CI_0:= C_b(\overline{U_\CI})$ and $\CJ_0:=C_b(\pr(A) \setminus \overline{U_\CI})$ as disjoint ideals of $C_b(\pr(A))$ with their sum equals $C_b(\pr(A))$.
Thus, $\CJ_0 = \CI_0^\bot$ and $\CI_0 = \CJ_0^\bot$ (here $\bot$ are taken in $C_b(\pr(A))$).
Since $\CI\subseteq \CI_0$ and $\CI^\bot \subseteq \CJ_0$, we see that $\CJ_0 = \CI^\bot$ as well as $\CI_0 = \CI^{\bot\bot}$.
Thus, $ZM(A) = \CI^\bot + \CI^{\bot\bot}$ and the spectrum of $ZM(A)$ is extremely disconnected (by the equivalence of Statements (1) and (2), which was established above, for $ZM(A)$).

Let $J$ be a normal ideal of $M(A)$.
By Statement (4') (which was shown to be equivalent to Statement (1) in the above), there is $q\in \proj(ZM(A))\cup \{0\}$ with $J = qM(A)$, and hence
$$J\cap ZM(A) = qZM(A),$$
which is a normal ideal of $ZM(A)$.
The injectivity of $\Delta$ follows from the fact that if $p,q\in \proj(ZM(A))\cup \{0\}$ with $pZM(A) = qZM(A)$, then $p=q$.

On the other hand, we claim that $\Delta$ is surjective when $ZM(A)$ is an $AW^*$-algebra.
In fact, since $\pr(ZM(A))$ is extremely disconnected, we know from Statement (4) (which was shown to be equivalent to Statement (1)) for $ZM(A)$ that any normal ideal $J_0$ of $ZM(A)$ comes from a projection $q_0$ in $ZM(A)$, and we have $\Delta(q_0M(A)) = J_0$.

\smallskip\noindent
$6)\Rightarrow 1)$.
By the equivalences established above, it suffices to verify Statement (4').
Suppose that $J$ is a non-zero normal ideal of $M(A)$.
By the assumption, $J\cap ZM(A)$ is a normal ideal and there is $q\in \proj(ZM(A))\cup \{0\}$ with $J \cap ZM(A) = qZM(A)$, because $ZM(A)$ is a commutative $AW^*$-algebra, and we can apply Statement (4)  for $ZM(A)$ (which was shown to be equivalent to $\pr(ZM(A))$ being extremely disconnected).
Since 
$qM(A) \cap ZM(A) = qZM(A) = J\cap ZM(A)$, the injectivity of $\Delta$ tells us that $J = qM(A)$.

\smallskip\noindent
$1)\Leftrightarrow 7)$.
This is precisely \cite[Proposition 2.9]{AM94}.
\end{prf}

\medskip

If $\pr(A)$ is extremely disconnected, we know from the above that there are bijection correspondences between $\proj(ZM(A))\cup \{0\}$, the set of normal ideals of $A$, the set of normal ideals of $M(A)$ as well as the set of normal ideals of $ZM(A)$.

\medskip

In the case when $A$ unital, the implication $1)\Rightarrow 6)$ also follows implicitly from \cite[Theorem 2.1]{Some96} and \cite[Proposition 3.2]{AS90}, but it will be easier to prove this directly than to recall the definition of ``quasi-standard'' $C^*$-algebras and their properties.

\medskip

Since the canonical map from the spectrum of a $C^*$-algebra $A$ to $\pr(A)$ induces a bijection from the collection of open subsets (respectively, closed subsets) of the spectrum to the collection of open subsets (respectively, closed subsets) of $\pr(A)$, we know that $\pr(A)$ is extremely disconnected if and only if $A$ has an extremely disconnected spectrum.

\medskip

By \cite[Example 2.10]{AM94}, all prime $C^*$-algebras, all $AW^*$-algebras and local multipliers algebras of all $C^*$-algebras will satisfy the equivalent conditions in Proposition \ref{prop:cent-AW-st}.
Moreover, Proposition \ref{prop:cent-AW-st} also tells us that if $A$ has an extremely disconnected spectrum, then $A$ is prime if and only if $\dim ZM(A) = 1$.

\medskip

In the following, we denote by $\beta \KI$ the Stone-Cech compactification of a set $\KI$ (as a discrete topological space).

\medskip

\begin{prop}\label{prop:quasi-Ston}
Let $A$ be a $C^*$-algebra with an extremely disconnected spectrum.

\smnoind
(a) $A = A_{\dc,1} \oplus A_{\dc,\infty} \oplus A_{\tw,1} \oplus A_{\tw,\infty} \oplus A_\ttr$.

\smnoind
(b) Suppose that there is a set $\KI$ with $ZM(A)\cong \ell^\infty(\KI)$.
Then $\KI$ is the disjoint union of subsets $\KI_{\dc,1}$, $\KI_{\dc,\infty}$, $\KI_{\tw,1}$, $\KI_{\tw,\infty}$ and $\KI_\ttr$ such that for $\#$ equals $\dc,1$ or $\dc,\infty$ or $\tw,1$ or $\tw,\infty$ or $\ttr$, one can find a  continuous field $\{A_\omega\}_{\omega\in \beta\KI_\#}$ of $C^*$-algebras on $\beta\KI_\#$ with $A_\#$ being the algebra of continuous sections of $\{A_\omega\}_{\omega\in \beta\KI_\#}$.
Moreover, for each $i\in \KI_\# \subseteq \beta \KI_\#$, the $C^*$-algebra $A_i$ is prime and it is an ideal of $A_\#$. 
The canonical map $\Theta: A_\#\to\prod_{i\in \KI_\#}^{\ell^\infty} A_i$ is injective and $\overline{\sum_{i\in \KI_\#} A_i}$ can be regarded as an essential ideal of $A_\#$.
\end{prop}
\begin{prf}
(a) By Relation \ref{eqt:her(e)-bot}, Proposition \ref{prop:cent-AW-st}, as well as Theorems \ref{thm:inf-part} and \ref{thm:decomp}, the central open projections corresponding to the ideals $A_{\dc,1}$, $A_{\dc,\infty}$, $A_{\tw,1}$, $A_{\tw,\infty}$ and $A_\ttr$ belong to $M(A)$ (since these projections are both open and closed) with their sum being $1$ (because this sum belongs to $M(A)$ and is dense in $1$ by Corollary \ref{cor:decomp}).
This gives the conclusion.

\smnoind
(b) Let $p_\#$ be the projection in $ZM(A)$ with $A_\# = p_\# A$ (see Statement (4)).
Then $p_\#\in \ell^\infty(\KI)$ is the characteristic function for a subset $\KI_\#$, and $\KI$ is the disjoint union of  $\KI_{\dc,1}$, $\KI_{\dc,\infty}$, $\KI_{\tw,1}$, $\KI_{\tw,\infty}$ and $\KI_\ttr$ (by part (a)).
As $A =  A_{\dc,1} \oplus A_{\dc,\infty} \oplus A_{\tw,1} \oplus A_{\tw,\infty} \oplus A_\ttr$, we know that 
$$ZM(A) =  ZM(A_{\dc,1}) \oplus ZM(A_{\dc,\infty}) \oplus ZM(A_{\tw,1}) \oplus ZM(A_{\tw,\infty}) \oplus ZM(A_\ttr),$$
and $ZM(A_\#) = \ell^\infty(\KI_\#)$. 
Moreover, one may regard $A_\#$ as a Banach $\ell^\infty(\KI_\#)$-normed module, and hence  consists of the continuous sections of a continuous field $\{A_\omega\}_{\omega\in \beta\KI_\#}$ of $C^*$-algebras (see e.g.\ \cite{DG83}).

For every $i\in \KI_\#$, we consider $p_i$ to be the projection in $ZM(A)$ corresponding to the point mass at $i$.
Then we have 
$$A_i\cong A/(1-p_i)A \cong p_i A = \her_A(p_i).$$
Suppose that $J$ is a non-zero closed ideal of $A_i$ and $p\in \op(A)$ with $J = \her_A(p)$.
As $p\leq p_i$ and $p_i\in M(A)$ (which implies that $p_i$ is a closed projection of $A$), we know that $\overline{p}^1 \leq p_i$.
Since $\overline{p}^1\in M(A)$ (by Statement (3)) and $p_i$ is a minimal projection in $ZM(A)$, we know that $\overline{p}^1 = p_i$ and $J$ is an essential ideal of $A_i$.
This shows that $A_i$ is prime.

Clearly, $A_i$ is an ideal of $A_\#$ and $I_\#:=\overline{\sum_{i\in \KI_\#} A_i}$ is a closed ideal of $A_\#$.
Note that $A_\#$ also has an extremely disconnected spectrum.
There exists, by Statement (4), a projection $q\in ZM(A_\#)$ with $I_\#^\bot= qA_\#$. 
Hence, $q$ is the characteristic function in $\ell^\infty(\KI_\#)$ for a subset $S\subseteq \KI_\#$.
If $S$ is non-empty, then for each $i\in S$, we have the contradiction that $A_i\subseteq I_\#\cap I_\#^\bot$.
Thus, $I_\#$ is an essential ideal of $A_\#$.

Finally, suppose that $a\in \ker \Theta$.
If we consider $a$ as a continuous section on $\{A_\omega\}_{\omega\in \beta \KI_\#}$, then its values at every $i\in \KI_\#$ is zero and hence $a$ is the zero section (as $\KI_\#$ is dense in $\beta\KI_\#$).
\end{prf}

\medskip

Proposition \ref{prop:quasi-Ston}(a) tells us that the decomposition in Corollary \ref{cor:decomp} extends the corresponding one for $AW^*$-algebras.
Moreover, part (b) applies to all $AW^*$-algebras whose centers are of the form $\ell^\infty(\KI)$.

\medskip

\subsection{Standard $C^*$-algebras}
We end this paper by considering some situations under which the $C^*$-algebra can be ``decomposed'' as a continuous field of prime $C^*$-algebras.
Let us begin with the following result.

\medskip

\begin{prop}\label{prop:ess-simp-decomp}
If $A$ is discrete (respectively, essentially infinite or tracially infinite; see Definition \ref{prop:pure-inf>III}), there is an open dense subset $\Xi_\dc$ (respectively, $\Xi_\infty$ or $\Xi_\tsinf$) of $\pr(A)$ such that $A/P$ is discrete (respectively, essentially infinite or tracially infinite) for every $P$ in $\Xi_\dc$ (respectively, $\Xi_\infty$ or $\Xi_\tsinf$).
\end{prop}
\begin{prf}
Suppose that $A$ is discrete.
By \cite[Theorem 2.3(v)]{PZ00}, the largest type $\tI$ ideal $A_{\tI}$ is essential.
Thus,
$$\Xi_\dc:=\pr(A)\setminus \hull(A_{\tI})$$
is an open dense subset of $\pr(A)$.
For every $P\in \Xi_\dc$, the image of $A_{\tI}$ in $A/P$
is a non-zero closed ideal of type $\tI$.
Thus, by Proposition \ref{prop:prime>5types}, we know that the prime $C^*$-algebra $A/P$ is discrete.

Secondly, we consider the case when $A$ is essentially infinite.
By Proposition \ref{prop:equiv-anti-fin}(a), there is an essential closed ideal $J_0$ of $A$ with $T_s(J_0) = \emptyset$.
Set
$$\Xi_\infty:= \pr(A)\setminus {\hull}(J_0).$$
Then $\Xi_\infty$ is an open dense subset of $\pr(A)$.
If $P\in \Xi_\infty$ and $q:A\to A/P$ is the quotient map, then $q(J_0)\neq (0)$ and does not have a tracial state.
Thus, $q(J_0)$ is essentially infinite and so is the prime $C^*$-algebra $A/P$ (by Proposition \ref{prop:prime>5types}).

In the case when $A$ is tracially infinite, there is an essential closed ideal $J_1\subseteq A$ with $\T(J_1) = \{0\}$.
The set
$$\Xi_\tsinf:= \pr(A)\setminus {\hull}(J_1).$$
is open and dense in $\pr(A)$.
Let $P\in \Xi_\tsinf$ and consider $q:A\to A/P$ to the quotient map.
By the construction, $q(J_1)\neq (0)$.

Suppose that $\T(q(J_1))\neq \{0\}$.
By the argument of Proposition \ref{prop:pure-inf>III}(a), one can find $y\in q(J_1)_+$ with $T_s\big(\overline{yq(J_1)y}\big) \neq \emptyset$.
Thus, the hereditary $C^*$-subalgebra $q^{-1}\big(\overline{yq(J_1)y}\big)\cap J_1$ of $J_1$ (which is non-zero because $J_1$ is an essential ideal of $A$) has a tracial state $\tau$.
As in the proof of Proposition \ref{prop:pure-inf>III}(a), the lower semi-continuous trace extension of $\tau$ to $J_1$ (as given in \cite[Lemma 4.6]{CP79}) is semi-finite.
This gives a contradiction.

Now, since $\T(q(J_1)) = \{0\}$ and $q(J_1)$ is an essential ideal of $A/P$ (as $A/P$ is prime), we know that $A/P$ is tracially infinite.
\end{prf}

\medskip

\begin{eg}
Let $A=\CL(\ell^2)$ and $P= \CK(\ell^2)$.
Then $A$ is discrete and $A/P$ is of type $\ttr$.
Therefore, one cannot expect $\Xi_\dc = \pr(A)$ in Proposition \ref{prop:ess-simp-decomp}(b).
\end{eg}

\medskip

By the Dauns-Hofmann theorem (see e.g. \cite[Corollary 4.4.8]{Ped79}), one has a continuous map $\Phi: \pr(A) \to \pr(ZM(A))$ with dense range.
In fact, if we identify elements in $\pr(ZM(A))$ with one-dimensional unital $^*$-representations (i.e., characters) of $ZM(A)$, then $\Phi(P)(x) = \|x+\ti P\|_{M(A)/\ti P}$ for every $x\in ZM(A)_+$, where is $\ti P$ the unique element in $\pr(M(A))$ with $\ti P\cap A = P$, i.e., 
$$\ti P := \{m\in M(A): mA, Am \subseteq P\}.$$
Thus, if $\omega$ is a character on $ZM(A)$, then $\Phi(P) = \omega$ if and only if $\ker(\omega)\subseteq \ti P$.

\medskip

Let $I$ be a closed ideal of $A$.
It is well-known that $I$ is the algebra of bounded continuous sections of a continuous fields $\Upsilon^I$ of $C^*$-algebras over $\pr(ZM(A))$ with the fiber $\Upsilon^I_\omega$ over $\omega\in \pr(ZM(A))$ being $I/\overline{I\ker(\omega)}$ (see e.g.\ \cite{DG83}).

\medskip

The following lemma gives another description for the fiber and it also gives a ``visual description'' of the image of an open subset of $\pr(A)$ under $\Phi$.
This lemma could be known, but since we do not find it in the literature, we give it argument here for completeness.

\medskip

\begin{lem}\label{lem:image-open-dense}
Let $I\subseteq A$ be a non-zero closed ideal. 
Consider any $\omega_0\in \pr(ZM(A))$ and $P_0\in \pr(A)\setminus \hull(I)$. 

\smnoind
(a) $\Phi(P_0) = \omega_0$ if and only if
$I\ker(\omega_0) \subseteq P_0$.

\smnoind
(b) $\overline{I\ker(\omega_0)} = \bigcap\{P\in \pr(A)\setminus \hull(I): \Phi(P) = \omega_0\}$. 

\smnoind
(c) $\Phi\big(\pr(A)\setminus \hull(I)\big) = \big\{\omega\in \pr(ZM(A)): \Upsilon^I_\omega \neq (0)\big\}$. 
\end{lem}
\begin{prf}
(a) If $\omega_0 = \Phi(P_0)$, then $\ker(\omega_0)\subseteq \ti P_0$, and we have  $I\ker(\omega_0)\subseteq P_0$.
Conversely, suppose that $I\ker(\omega_0)\subseteq P_0$.
Let $Q_0:= P_0\cap I\in \pr(I)$ and
$$\ti Q_0:= \{m\in M(I): mI,Im\subseteq Q_0\}.$$
Consider $\Psi: M(A)\to M(I)$ to be the canonical $^*$-homomorphism.
Set $J:= A\Psi^{-1}\big(\ti Q_0\big) \subseteq A$.
As $I\ker(\omega_0)\subseteq P_0\cap I$, one has $\Psi(\ker(\omega_0))\subseteq \ti Q_0$, and so 
$$A\ker(\omega_0)\subseteq J.$$
Moreover, it follows from $\Psi(IJ)\subseteq I \ti Q_0 \subseteq Q_0$ that $IJ \subseteq Q_0\subseteq P_0$ (since $\Psi$ restricts to the identity map  on $I$).
This implies that $J\subseteq P_0$ (as $I\nsubseteq P_0$ by the assumption).
Thus, $A\ker(\omega_0) \subseteq P_0$, which implies that $\ker(\omega_0)\subseteq \ti P_0$ and so, $\omega_0 = \Phi(P_0)$.

\smnoind
(b) Since $\overline{I\ker(\omega_0)}$ is a closed ideal of $I$, one has $$\overline{I\ker(\omega_0)} = \bigcap \{P\in \pr(A)\setminus \hull(I): \overline{I\ker(\omega_0)}\subseteq P\},$$
and this part follows from part (a).

\smnoind
(c) It follows from part (a) that $\omega \in \Phi\big(\pr(A)\setminus \hull(I)\big)$ if and only if there exists $P\in \pr(A)\setminus \hull(I)= \pr(I)$ such that $I\ker(\omega) \subseteq P$, which is equivalent to the closed ideal $\overline{I\ker(\omega)}$ of $I$ being proper (in other words, $I/\overline{I\ker(\omega)} \neq (0)$).
\end{prf}

\medskip

Following \cite[p.351]{AS90} and \cite[p.125]{Some96}, for each $\omega\in \pr(ZM(A))$, we call
$$G(\omega):=\bigcap \Phi^{-1}(\omega)$$
a \emph{Glimm ideal} of $A$ (here we use the usual convention that $\bigcap \emptyset = A$).

\medskip

\begin{thm}\label{thm:Glimm>prime}
Suppose that the map $\Phi$ is open and every non-zero Glimm ideal of a $C^*$-algebra $A$ is a prime ideal.

\smnoind
(a) $A$ is the algebra of $C_0$-sections of a continuous field of $C^*$-algebras over $\Omega_A:=\Phi(\pr(A))$ with each fiber being non-zero and of one of the five types.

\smnoind
(b) If $A$ is discrete (respectively, essentially infinite or tracially infinite), there is an open dense subset $\Omega_\dc$ (respectively, $\Omega_\infty$ or $\Omega_\tsinf$) of $\Omega_A$ such that the fiber over elements in $\Omega_\dc$ (respectively, $\Omega_\infty$ or $\Omega_\tsinf$) are discrete (respectively, essentially infinite or tracially infinite).
\end{thm}
\begin{prf}
(a) Consider any $\omega\in \Omega_A$. 
By Lemma \ref{lem:image-open-dense}(b), one has 
$$\overline{A\ker(\omega)} = G(\omega),$$ 
and hence is a prime ideal by the hypothesis.
Moreover, it follows from Lemma \ref{lem:image-open-dense}(c) that $A/\overline{A\ker(\omega)} = \Upsilon^A_\omega\neq (0)$. 
Now, Proposition \ref{prop:prime>5types} tells us that they are of one of the five types.

\smnoind
(b) We will only consider the case when $A$ is discrete since the other two cases follow from similar arguments. 
Suppose that $A_{\tI}$ is the largest type $\tI$ closed ideal of $A$ and $\Xi_\dc:= \pr(A)\setminus {\hull}(A_{\tI})$.
By the proof of Proposition \ref{prop:ess-simp-decomp}, $\Xi_\dc$ is an open dense subset of $\pr(A)$.
Set $\Omega_\dc:=\Phi(\Xi_\dc)$.
Then $\Omega_\dc$ is an open dense subset of $\Omega_A$ (by the hypothesis) satisfying
\begin{equation*}\label{eqt:Omega-d}
\Omega_\dc\subseteq \{\omega\in \Omega_A: G(\omega)\subseteq P, \text{ for some } P\in \Xi_\dc\}
\end{equation*}
(since $G(\omega)\subseteq P$ whenever $P\in \Phi^{-1}(\omega)$).

Pick any $\omega\in \Omega_d$.
If  $A_{\tI}$ is contained in $G(\omega)$, then we know that $A_{\tI}\subseteq P$ for some $P\in \pr(A)\setminus {\hull}(A_{\tI})$, which is absurd.
Consequently, for every $\omega\in \Omega_d$, the image of $A_{\tI}$ in the prime $C^*$-algebra $\Upsilon^A_\omega = A/G(\omega)$ is a non-zero type $\tI$ closed ideal, which implies that $\Upsilon^A_\omega$ is discrete (because of Proposition \ref{prop:prime>5types}).
\end{prf}

\medskip

If $\pr(A)$ is Hausdorff, then $\Phi$ is injective and the hypothesis of Theorem \ref{thm:Glimm>prime} is satisfied.
Moreover, if $A$ is a ``standard'' $C^*$-algebras as in \cite[p.127]{Some96}, then the two hypothesis in Theorem \ref{thm:Glimm>prime} are satisfied (see also \cite[Theorem 3.4]{AS90}). %changed after submission
Thus, Theorem \ref{thm:Glimm>prime} gives the following corollary, thanks to \cite[p.127]{Some96} and \cite[Theorem 2.4]{Some96}.

\medskip

\begin{cor}\label{cor:ess-simp}
If either $\pr(A)$ is Hausdorff, or $A$ is an $AW^*$-algebra, or $A$ is the local multiplier algebra of another $C^*$-algebra, then $A$ is the algebra of continuous $C_0$-sections of a continuous fields of prime $C^*$-algebras over a locally compact Hausdorff space $\Omega_A$ with each fiber being non-zero and of one of the five types.
If, in addition, $A$ is discrete (respectively, essentially infinite or tracially infinite), there is an open dense subset of $\Omega_A$ on which the fiber are
nonzero and discrete (respectively, essentially infinite or tracially infinite).
\end{cor}

\medskip

\bibliographystyle{plain}

\begin{thebibliography}{99}
\bibitem{Akemann69}
C.\ A.\ Akemann, The general Stone-Weierstrass problem, J.\ Funct.\ Anal., \textbf{4} (1969), 277-294.

\bibitem{Akemann70}
C.\ A.\ Akemann, Left ideal structure of $C^*$-algebras, J.\ Funct.\ Anal., \textbf{6} (1970), 305-317.

\bibitem{AM94}
P.\ Ara and M.\ Mathieu, On the central Haagerup tensor product, Proc.\ Edin.\ Math.\ Soc., \textbf{37} (1994), 161-174.

\bibitem{AS90}
R.\ J.\ Archbold and D.\ W.\ B.\ Somerset, Quasi-standard C*-algebras, Math.\ Proc.\ Comb.\ Phil.\ Soc., \textbf{107} (1990), 349-360.

\bibitem{Blac}
B.\ Blackadar, Operator algebras - theory of $C^*$-algebras and $W^*$-algebras, in \emph{Operator Algebras and Non-commutative Geometry III},
Encyc.\ Math.\ Sci., \textbf{122}, Springer-Verlag, Berlin (2006).

\bibitem{BD}
F.F. Bonsall and J. Duncan, \emph{Complete normed algebras}, Springer-Verlag, Berlin (1973).

\bibitem{BP91}
L.\ G.\ Brown and G.\ K.\ Pedersen,
$C^*$-algebras of real rank zero, J.\ Funct.\ Anal., \textbf{99} (1991), 131-149.

\bibitem{Cu81}
J. Cuntz, $K$-theory for certain $C^*$-algebras, Ann. of Math. \textbf{113} (1981), 181-197.

\bibitem{CP79}
J.\ Cuntz and G.\ K.\ Pedersen,
Equivalence and traces on $C^*$-algebras, J.\ Funct.\ Anal., \textbf{33} (1979), 135-164.

\bibitem{Dix77}
J. Dixmier, \emph{$C^*$-Algebras}, North-Holland Publishing (1977).

\bibitem{DG83}
M.\ J.\ Dupr\'{e} and R.\ M.\ Gillette, \emph{Banach bundles, Banach
modules and automorphisms of $C\sp{*} $-algebras}, Research Notes in
Mathematics 92, Pitman (1983).

%\bibitem{Fell61}
%J.\ M.\ G.\ Fell, The structure of algebras of operator fields, Acta Math.\, \textbf{106} (1961), 233-280.

\bibitem{KRII}
R.\ V.\ Kadison and J.\ R.\ Ringrose, \emph{Fundamentals of the theory of operator algebras Vol. II: Advanced theory}, Pure and Applied Mathematics \textbf{100}, Academic Press (1986).

\bibitem{KR00}
E.\ Kirchberg and M.\ R{\o}rdam, Non-simple purely
infinite $C^*$-algebras, Amer.\ J.\ Math., \textbf{122} (2000), 637--666.

\bibitem{Lan}
E.\ C.\ Lance, \emph{Hilbert $C^*$-modules - A toolkit for operator algebraists}, Lond.\ Math.\ Soc.\ Lect.\ Note Ser.\ \textbf{210},  Camb.\ Univ.\ Press (1995).

\bibitem{LN}
C. W. Leung and C. K. Ng, Invariant ideals of twisted
crossed products, Math. Zeit. \textbf{243} (2003), 409-421.

\bibitem{LinHX90}
H.\ Lin, Equivalent open projections and corresponding hereditary
$C^*$-subalgebras, J.\ Lond.\ Math.\ Soc., \textbf{41}
(1990), 295-301.

\bibitem{Murphy}
G.J. Murphy, \emph{$C^*$-algebras and operator theory}, Academic Press (1990).

\bibitem{Murray90}
F.\ J.\ Murray,
The rings of operators papers,
in \emph{The legacy of John von Neumann} (Hempstead, NY, 1988), 57-60, Proc.\ Sympos.\ Pure Math.\ \textbf{50}, Amer.\ Math.\ Soc., Providence, R.I.\ (1990).

\bibitem{Murray-vonNeumann36}
F.\ J.\ Murray and J.\ Von Neumann,
On rings of operators,
Ann.\ of Math.\ (2), \textbf{37} (1936), 116-229.

\bibitem{Ng-str-amen-rep}
C.\ K.\ Ng, Strictly amenable representations of reduced group $C^*$-algebras, Int. Math. Res. Notices, to appear.

\bibitem{Ng-ext-lscsf-tr}
C.\ K.\ Ng, On extension of lower semi-continuous semi-finite traces, in preparation.

\bibitem{NW-Mv-class}
C.\ K.\ Ng and N.\ C.\ Wong, A Murray-von Neumann type classification of $C^*$-algebras, in \emph{Operator Semigroups Meet Complex Analysis, Harmonic Analysis and Mathematical Physics, Herrnhut, Germany (in honor of Prof. Charles Batty for his 60th birthday)}, Operator Theory: Advance and Applications, \textbf{250}, Springer Internat. Publ. (2015), 369-395.

%\bibitem{NW-ext-v-alg-class}
%C. K. Ng and N. C. Wong, An extension of $W^*$-algebra classification to $C^*$-algebras, preprint.

%\bibitem{ORT}
%E.\ Ortega, M.\ R{\o}rdam and H.\ Thiel,
%The Cuntz semigroup and comparison of open projections, J.\ Funct.\ Anal., \textbf{260} (2011), 3474--3493.

\bibitem{Ped79}
G.\ K.\ Pedersen,
\emph{$C^*$-algebras and their automorphism groups},
Academic Press (1979).

\bibitem{PZ00}
C.\ Peligrad and L.\ Zsid\'{o}, Open projections of $C^*$-algebras:
Comparison and Regularity, in \emph{Operator Theoretical Methods, 17th
Int.\ Conf.\ on Operator Theory, Timisoara (Romania), June 23-26,
1998}, Theta Found.\ Bucharest (2000), 285-300

\bibitem{Rord03}
M. R\o rdam, A simple $C^*$-algebra with a finite and an infinite projection, Acta Math. \textbf{191} (2003), 109-142.

\bibitem{Some96}
D.\ W.\ B.\ Somerset, The local multiplier algebra of a C*-algebra, Quart.\ J.\ Math.\ Oxford, Ser.\ 2, \textbf{47} (1996), 123-132.

\bibitem{Zhang89}
S.\ Zhang, Stable isomorphism of
hereditary $C^\ast$-subalgebras and stable equivalence of open projections,
Proc.\ Amer.\ Math.\  Soc., \textbf{105} (1989), 677-682.
\end{thebibliography}

\end{document}